\documentclass[12pt]{article}

\usepackage{
	enumitem,
	times,
	amsmath,
	amssymb,
	amsthm,
	mathtools,
	fancyhdr,
	color,
	thmtools,
	etoolbox,
	xspace,
	hyperref,
	lscape,
	thm-restate
}

\usepackage[margin=2cm]{geometry}

\declaretheorem[style=definition,numberwithin=section]{definition}
\declaretheorem[style=definition,qed=$\oslash$,sibling=definition]{example}

\declaretheorem[style=slanted,sibling=definition]{theorem}
\declaretheorem[style=slanted,sibling=definition]{conjecture}
\declaretheorem[style=slanted,sibling=definition]{lemma}
\declaretheorem[style=slanted,sibling=definition]{proposition}

\declaretheorem[style=definition]{claim}
\declaretheorem[style=definition,sibling=example]{remark}

\AtEndEnvironment{proof}{\setcounter{claim}{0}}

\usepackage[square, comma, numbers, sort&compress ]{natbib}
\newcommand{\C}{\ensuremath{\mathbb{C} }}
\newcommand{\R}{\ensuremath{\mathbb{R} }}
\newcommand{\N}{\ensuremath{\mathbb{N} }}

\DeclareMathOperator{\Tr}{Tr}
\DeclareMathOperator{\End}{End}
\DeclareMathOperator{\Ad}{Ad}

\DeclareMathOperator{\Span}{span}

\newcommand\scalemath[2]{\scalebox{#1}{\mbox{\ensuremath{\displaystyle #2}}}}
\title{On the Mathieu Conjecture for $SU(N)$ and $SO(N)$}
\author{Kevin Zwart\footnote{Email address: \href{mailto:kevin.zwart@ru.nl}{kevin.zwart@ru.nl}}}
\date{%
	Radboud University Nijmegen, IMAPP \\[2ex]%
	\today}

\begin{document}

\maketitle

\abstract{Building on work of M.\ M\"uger and L.\ Tuset, we reduce the Mathieu
conjecture, formulated by O.\ Mathieu in 1997, for $SU(N)$ to a simpler conjecture in purely abelian
terms. We sketch a similar reduction for $SO(N)$. The proofs rely on
Euler-style parametrizations of these groups, which we discuss
including proofs.}

\section{Introduction}

One of the famous conjectures is the Jacobi conjecture:
\begin{conjecture}[Jacobi Conjecture]
	Let $f:\C^n\rightarrow\C^n$ be a polynomial map, which is to say that the component $f_i:\C^n\rightarrow\C$ is a polynomial map for each $i=1,\ldots,n$. If the derivative $D_xf$ is invertible for all $x\in \C^n$, then $f$ is bijective and its inverse is a polynomial map as well.
\end{conjecture}
It was conjectured by Keller in 1939, and although many attempts have been made, no proofs have been found yet (for more details about the Jacobi conjecture and its history, see e.g. \cite{smale,bass}). In 1997, O. Mathieu published a paper \cite{Mathieu} in which he proved the Jacobi conjecture if another conjecture was to hold. He considered a connected compact Lie group $G$, and states that for any finite-type function $f,h$ on $G$, if the integral $\int_G f(g)^P dg=0$ for all $P\in \N$, then $\int_G f(g)^Ph(g) dg=0$ for large enough $P$.

In 1998, only one year after the publication of Mathieu's paper, Duistermaat and van der Kallen \cite{Duistermaat} published an article where they proved Mathieu's conjecture in the case of an abelian connected compact group. Although no proofs are known for the non-abelian case, some attempts were made. A paper by Dings and Koelink \cite{Dings-Koelink} tried to prove the Mathieu conjecture by expressing the finite-type functions by explicit matrix coefficients. Influenced heavily by this paper, M\"uger and Tuset published a paper \cite{Mueger} which reduces the Mathieu conjecture on $SU(2)$ to a conjecture which involves only Laurent polynomials on $\C$ with coefficients that are polynomial functions on $\R$ itself.

The goal of the present paper is to generalize the paper by M\"uger and Tuset to the compact matrix groups $SU(N)$ and $SO(N)$, where $N\geq 2$. A key ingredient to achieving this will be a generalization of the Euler decomposition. The Euler decomposition on $SU(2)$ has been known for some time but is mostly used by physicists under the name of Euler angles. This is no different in the case of $SU(N)$. Several (similar but not equal) versions of the Euler decomposition for $SU(N)$ exist, see for example Bertini et al. \cite{bertini2006euler}, Cacciatori et al. \cite{Cacciatori}, or Tilma et al. \cite{Parametrization_SU(N)}. In a similar way there exist several decompositions of $SO(N)$, see for example \cite{SO(N)_parametrization_raffenetti1970,SO(N)_hoffman1972,SO(N)_WIGNER}.

In our paper, we will reduce the Mathieu conjecture to a conjecture similar to that of M\"uger and Tuset \cite{Mueger}. We start by looking at the matrix coefficients of the generalized Euler decomposition on $SU(N)$ and $SO(N)$, and we find that any finite-type function can be described by a function on $\C^n\times \R^k$. To be more specific, any finite type function reduces to a function $f:\C^n\times \R^k\rightarrow\C$ which can be written as $f(z,x)=\sum_{\vec{m}}c_{\vec{m}}(x)z^{\vec{m}}$ where $\vec{m}=(m_1,\ldots,m_k)$ is a multi-index where $m_i\in\bigcup_{j=1}^N\frac{1}{j}\mathbb{Z}$ for each $i$, and $c_{\vec{m}}(x)$ is a polynomial in $x_1,\ldots,x_k$ and $\sqrt{1-x_1^2},\ldots,\sqrt{1-x_k^2}$. Assuming these functions satisfy other conjectures, the Mathieu conjecture is proven for $SU(N)$ and $SO(N)$. The proof uses the explicit description of the Euler decomposition on $SU(N)$ and $SO(N)$ and the properties of the Haar measure in these parametrizations. In Section \ref{sec:SU(N)} we will focus on the group $SU(N)$, while in Section \ref{sec:SO(N)} the group $SO(N)$ will be considered. The final part of the paper is dedicated to proving the generalized Euler decomposition we used throughout this paper, with the corresponding explicit description of the Haar measure in this parametrization.

\textit{Acknowledgment:} The author would like to thank Michael M\"uger for proposing the subject, and the many valuable discussions we had. He also wishes to thank Erik Koelink for feedback and suggestions.

\section{The case of $SU(N)$}\label{sec:SU(N)}
In this paper we will reduce the Mathieu's conjecture on $SU(N)$ and $SO(N)$ with $N \geq 2$. We start by recalling Mathieu's conjecture. To do so, we first introduce the notion of a finite-type function:
\begin{definition}
	Let $G$ be a compact Lie group. A function $f:G\rightarrow\C$ is called a \emph{finite-type function} if it can be written as a finite sum of matrix components of irreducible continuous representations, i.e. $$f(x)=\sum_{j=1}^N\Tr\left[a_j\pi_j(x)\right],$$ where $(\pi_j,V_j)$ is an irreducible continuous representation of $G$, and  $a_j\in\End(V_j)$.  
\end{definition}

\begin{conjecture}[The Mathieu Conjecture \cite{Mathieu}]
	Let $G$ be a compact connected Lie group. If $f,h$ are finite-type functions such that $\int_G f^P \,dg=0$ for all $P\in\N$, then $\int_G f^Ph \,dg=0$ for all large enough $P$.
\end{conjecture}

In this section we will focus on $SU(N)$. We will base our parametrization and Haar measure on \cite{Parametrization_SU(N),Parametrization_SU(4)}. For completeness, we included an appendix dedicated to proving the parametrization.

For simplicity, we will define the generators of $\mathfrak{su}(N)$ for $N\in\N$. Let $j=1,2,\ldots,N-1$ and $k=1,2,\ldots,2j$ and define the matrices $\lambda_j\in \mathfrak{su}(N)$ in the following way\footnote{In most physics papers the matrices $\{i\lambda_j\}_j$ are called Gell-Mann matrices, see e.g. \cite{Parametrization_SU(N),Parametrization_SU(4),bertini2006euler}}
\begin{align*}
[\lambda_{j^2-1+k}]_{\mu,\nu}&:=i(\delta_{\lceil\frac{k}{2}\rceil,\mu}\delta_{j+1,\nu}+\delta_{j+1,\mu}\delta_{\lceil\frac{k}{2}\rceil,\nu})\quad\qquad\text{if $k$ is odd},\\	[\lambda_{j^2-1+k}]_{\mu,\nu}&:=\delta_{\frac{k}{2},\mu}\delta_{j+1,\nu}-\delta_{j+1,\mu}\delta_{\frac{k}{2},\nu}\qquad\qquad\quad\text{if $k$ is even},\\ 
[\lambda_{(j+1)^2-1}]_{\mu,\nu}&:=\begin{pmatrix}
i\mathbf{1}_j&0&&&\\
0&-ij&&&\\
&&0&&\\
&&&\ddots&\\
&&&&0
\end{pmatrix},
\end{align*}
where $\mathbf{1}_j$ is the $j\times j$ identity matrix. The matrices $\lambda_1\ldots,\lambda_{N^2-1}$ span $\mathfrak{su}(N)$. For example, the first eight matrices are given by 
\begin{align*}
\lambda_{1}=\begin{pmatrix}
0&i&0&\ldots&0\\
i&0&0&\ldots&0\\
0&0&0&\ldots&0\\
\vdots&\vdots&\vdots&\ddots&\vdots\\
0&0&0&\ldots&0
\end{pmatrix},\qquad \lambda_{2}=\begin{pmatrix}
0&1&0&\ldots&0\\
-1&0&0&\ldots&0\\
0&0&0&\ldots&0\\
\vdots&\vdots&\vdots&\ddots&\vdots\\
0&0&0&\ldots&0
\end{pmatrix},\qquad\lambda_3=\begin{pmatrix}
i&0&0&\ldots&0\\
0&-i&0&\ldots&0\\
0&0&0&\ldots&0\\
\vdots&\vdots&\vdots&\ddots&\vdots\\
0&0&0&\ldots&0
\end{pmatrix},\\
\lambda_{4}=\begin{pmatrix}
0&0&i&\ldots&0\\
0&0&0&\ldots&0\\
i&0&0&\ldots&0\\
\vdots&\vdots&\vdots&\ddots&\vdots\\
0&0&0&\ldots&0
\end{pmatrix},\qquad \lambda_{5}=\begin{pmatrix}
0&0&1&\ldots&0\\
0&0&0&\ldots&0\\
-1&0&0&\ldots&0\\
\vdots&\vdots&\vdots&\ddots&\vdots\\
0&0&0&\ldots&0
\end{pmatrix},\qquad\lambda_6=\begin{pmatrix}
0&0&0&\ldots&0\\
0&0&i&\ldots&0\\
0&i&0&\ldots&0\\
\vdots&\vdots&\vdots&\ddots&\vdots\\
0&0&0&\ldots&0
\end{pmatrix},\\
\lambda_{7}=\begin{pmatrix}
0&0&0&\ldots&0\\
0&0&1&\ldots&0\\
0&-1&0&\ldots&0\\
\vdots&\vdots&\vdots&\ddots&\vdots\\
0&0&0&\ldots&0
\end{pmatrix},\qquad\lambda_8=\begin{pmatrix}
i&0&0&\ldots&0\\
0&i&0&\ldots&0\\
0&0&-2i&\ldots&0\\
\vdots&\vdots&\vdots&\ddots&\vdots\\
0&0&0&\ldots&0
\end{pmatrix}.
\end{align*}
With this basis of $\mathfrak{su}(N)$ we have the following lemma:

\begin{restatable}[Generalized Euler Angles]{lemma}{EulerAngles}\label{lemma:Euler_Angles}
	Let $N\geq 2$. Define inductively the mapping $F_N:([0,\pi]\times[0,2\pi]^{N-2})\times([0,\pi]\times[0,2\pi]^{N-3})\times\cdots\times([0,\pi]\times [0,2\pi])\times[0,\pi]\times \left[0,\frac{\pi}{2}\right]^{\frac{N(N-1)}{2}}\times[0,2\pi]\times\cdots\times\left[0,\frac{2\pi}{N-1}\right]\rightarrow SU(N)$ by $F_1\equiv 1$ and
	\begin{equation}\label{eq:Euler_parametrization}
	\begin{split}
	&F_N(\phi_1,\ldots\phi_{\frac{N(N-1)}{2}},\psi_1,\ldots,\psi_{\frac{N(N-1)}{2}},\omega_1,\ldots,\omega_{N-1}):=\\
	&\left(\prod_{2\leq k\leq N}A(k)(\phi_{k-1},\psi_{k-1})\right)\cdot\begin{pmatrix}
	F_{N-1}(\phi_{N},\ldots,\phi_{\frac{N(N-1)}{2}},\psi_{N},\ldots,\psi_{\frac{N(N-1)}{2}},\omega_1,\ldots,\omega_{N-2})&0\\
	0&1
	\end{pmatrix} e^{\lambda_{N^2-1}\omega_{N-1}},
	\end{split}	
	\end{equation}
	where $A(k)(x,y):=e^{\lambda_{3}x}e^{\lambda_{(k-1)^2+1}y}$, and $\psi_j\in\left[0,\frac{\pi}{2}\right],\,\omega_j\in\left[0,\frac{2\pi}{j}\right]$ for all $j$. Here we denoted the product as $$\prod_{2\leq k\leq N}A(k)(\phi_{k-1},\psi_{k-1}):=A(2)(\phi_1,\psi_1)\cdot \cdots\cdot A(N)(\phi_{N-1},\psi_{N-1}).$$ This mapping is surjective. Moreover it is a diffeomorphism on the interior of the hypercube.
\end{restatable}

\begin{remark}\label{remark:KAK_decomposition}
	This lemma tells us that we have a parametrization of $SU(N)$ up to a measure zero sets. In the case of $SU(2)$, this reduces to the Euler angles parametrization, which is given by $$SU(2)=\left\{\left.\begin{pmatrix}
	e^{i\phi}&0\\
	0&e^{-i\phi}
	\end{pmatrix}\begin{pmatrix}
	\cos(\psi)&\sin(\psi)\\
	-\sin(\psi)&\cos(\psi)
	\end{pmatrix}\begin{pmatrix}
	e^{i\omega}&0\\
	0&e^{-i\omega}
	\end{pmatrix}\,\right|\, \phi\in[0,\pi],\psi\in[0,\pi/2], \omega\in [0,2\pi]\right\}.$$ To give further motivation for this parametrization, we can define a Cartan involution $\theta$ by $$\theta:\mathfrak{su}(N)\rightarrow\mathfrak{su}(N),\qquad\theta=\Ad(e^{\frac{\pi}{N-1}\lambda_{N^2-1}})=\Ad\left(\begin{pmatrix}
	e^{\frac{i\pi}{N-1}} \mathbf{1}_{N-1}&0\\
	0&-1
	\end{pmatrix}\right).$$ We see that $\theta=1$ on the subalgebra $\mathfrak{k}$ spanned by $\lambda_1,\ldots,\lambda_{(N-1)^2-1}$ and $\lambda_{N^2-1}$, and $\theta=-1$ on the vector space $\mathfrak{p}:=\Span_\R(\lambda_{(N-1)^2},\ldots,\lambda_{N^2-2})$. In addition, note that $\mathfrak{k}\simeq \mathfrak{su}(N-1)\oplus \mathfrak{u}(1)$. Since $SU(N)$ is connected for all $N$, the corresponding connected Lie group $K$ such that $\text{Lie}(K)=\mathfrak{k}$ can be seen as $$K=\begin{pmatrix}
	SU(N-1)&0\\
	0&1
	\end{pmatrix}e^{\R\lambda_{N^2-1}}\simeq S(U(N-1)\times U(1)).$$ We choose the maximal abelian subalgebra $\mathfrak{a}\subset \mathfrak{p}$ as $\mathfrak{a}=\R\lambda_{(N-1)^2+1}$. The $KAK$-decomposition \cite{KnappBeyond} then gives
	$$SU(N)=KAK\simeq [S(U(N-1)\times U(1))]e^{\R\lambda_{(N-1)^2+1}}[S(U(N-1)\times U(1))].$$
	Our lemma states that, up to a measure zero set, there exists a subset $L\subset K$ such that $SU(N)$ is diffeomorphic to $LAK$. We also note that by construction $SU(N)/K$ is a symmetric space and is diffeomorphic to the complex projective plane $\mathbb{CP}^{N-1}$.
\end{remark}

\begin{restatable}{lemma}{HaarMeasure}\label{lemma:Haar_Measure}
	Let $N\geq 2$ and $F_N$ be the Euler parametrization of $SU(N)$. The Haar measure $dg_{SU(N)}$ is then given inductively by $$dg_{SU(2)}=C_2\cos(\psi_1)\sin(\psi_1)\,d\phi_1d\psi_1d\omega_1$$ and
	\begin{align}\label{eq:Haar_Measure}
	\begin{split}dg_{SU(N)}=C_N\cos(\psi_{N-1})\sin^{2(N-1)-1}(\psi_{N-1})&\left[\prod_{j=1}^{N-2}\cos^{2j-1}(\psi_{j})\sin(\psi_j)\right]\cdot\\
	&d\phi_1\ldots d\phi_{N-1}d\psi_1\ldots d\psi_{N-1} dg_{SU(N-1)}d\omega_N,\end{split}
	\end{align}
	where $C_n:= \frac{(n-1)!(n-1)}{2\pi^n}$ for all integers $n\geq 2$.
\end{restatable}

As mentioned, Lemma \ref{lemma:Euler_Angles} and \ref{lemma:Haar_Measure} are proved in the appendix. With these lemmas, we can start the discussion of Mathieu's conjecture. Note that any finite-type function on $SU(N)$ is a sum of products of matrix coefficients since the irreducible representations of $SU(N)$ are polynomials in the matrix coefficients. By the parametrization given in Equation (\ref{eq:Euler_parametrization}), we see that these products consist of (powers of) $\sin(\psi_j), \cos(\psi_k),e^{i\phi_l}$ and $e^{i\omega_m}$. Therefore any finite-type function $h$  can be written as
\begin{align}\begin{split}\label{eq:K-finite_function}
h(g)= \sum_{j=1}^M\sum_{i=1}^Q c_{ij}& e^{ik_{ij}^1\phi_1}\sin^{m_{ij}^1}(\psi_1)\cos^{n_{ij}^1}(\psi_1)\cdots e^{ik_{ij}^{N-1}\phi_{N-1}}\sin^{m_{ij}^{N-1}}(\psi_{N-1})\cos^{n_{ij}^{N-1}}(\psi_{N-1})\\
&\cdot (h_{SU(N-1)})_{ij}(g_{SU(N-1)}) e^{il^N_{ij}\omega_{N-1}},\end{split}
\end{align}
where $g_{SU(N-1)}:=F_{N-1}(\phi_{N},\ldots,\phi_{\frac{N(N-1)}{2}},\psi_{N},\ldots,\psi_{\frac{N(N-1)}{2}},\omega_1,\ldots,\omega_{N-2})$ is the $SU(N-1)$ component of $g=F_N(\phi_1,\ldots,\omega_{N-1})$ as in Lemma \ref{lemma:Euler_Angles}, and $(h_{SU(N-1)})_{ij}$ is a family of finite-type functions on $SU(N-1)$. Also $k_{ij}^p,l_{ij}^p\in\mathbb{Z}$, $m_{ij}^p\in\N$ and $n_{ij}^p\in\{0,1\}$. We can achieve $n_{ij}^p\in \{0,1\}$ by using the equality $\cos^2(\psi_j)+\sin^2(\psi_j)=1$ repeatedly. Note that we sum over both $i$ and $j$. The sum over $i$ is to ensure we have all possible combinations of different terms, while the sum over $j$ allows for different powers of each term. For example, in $SU(2)$, we have the parametrization of the form $$SU(2)=\left\{e^{\phi_1\lambda_3}e^{\psi_1\lambda_1}e^{\omega_1\lambda_3}\,\Big|\,\phi_1\in [0,\pi],\psi_1\in \left[0,\frac{\pi}{2}\right],\omega_1\in[0,2\pi]\right\},$$ so any finite-type function is of the form
\begin{align*}
h_{SU(2)}(e^{\phi_1\lambda_3}e^{\psi_1\lambda_1}e^{\omega_1\lambda_3})=\sum_{j=1}^M c_je^{ik_j\phi_1}\sin^{m_j}(\psi_1)e^{il_j\omega_1}+c_j'e^{ik_j'\phi_1}\sin^{m_j'}(\psi_1)\cos(\psi_1)e^{il_j'\omega_1}.
\end{align*} 

\begin{remark}
	We note that if we restrict any finite type function $h$ to a closed subgroup $H$ of $SU(N)$, then $h|_{H}$ also is a finite-type function. This can easily be seen by the fact that any irreducible representation $(\pi,V)$ of $SU(N)$ is finite dimensional, hence $(\pi|_H,V)$ splits into finitely many irreducible representations $(\pi_{H,i}, V_i)$ of $H$, i.e. $\pi|_H\simeq \bigoplus_{i=1}^M \pi_{H,i}$. It is immediate then that $h|_H$ is again a finite-type function.
\end{remark}

\begin{lemma}\label{lemma:Main_Lemma}
	Let $h$ be a finite-type function on $SU(N)$ as in Equation (\ref{eq:K-finite_function}), and $N\geq 2$. Then for any $P\in \N$ we have
	\begin{align}
	\begin{split}\int_{SU(N)}h(g)^P\,dg=\frac{1}{2(N-1)i^{\frac{N(N+1)}{2}-1}}&\int_{[0,1]^\frac{N(N-1)}{2}}\int_{(S^*)^{\frac{N(N+1)}{2}-1}}\left[\widetilde{h_{SU(N)}}(x_1,\ldots,z_{\frac{N(N+1)}{2}-1})\right]^P\\
	&\cdot J_{SU(N)}(x_1,\ldots,x_{\frac{N(N-1)}{2}})\frac{dz_1}{z_1}\ldots\frac{dz_{\frac{N(N+1)}{2}-1}}{z_{\frac{N(N+1)}{2}-1}}dx_1\ldots dx_{\frac{N(N-1)}{2}}.\end{split}
	\end{align}
	Here $J_{SU(N)}$ is defined recursively by $J_{SU(1)}\equiv 1$ and, for $2\leq n\leq N$, by $$J_{SU(n)}(x_1,\ldots,x_{\frac{n(n-1)}{2}}):=C_n\,x_{n-1}^{2n-3}\left(\prod_{j=1}^{n-2}x_j(1-x_j^2)^{j-1}\right) J_{SU(n-1)}\left(x_{n},\ldots,x_{\frac{n(n-1)}{2}}\right),$$ where $C_n$ is as in Lemma \ref{lemma:Haar_Measure}, and where $\widetilde{h_{SU(N)}}$ is defined recursively by $\widetilde{h_{SU(1)}}\equiv1$ and by \begin{align}\label{eq:h_tilde}
	\widetilde{h_{SU(n)}}(x_1,\ldots,x_{\frac{n(n-1)}{2}},z_1,\ldots,&z_{\frac{n(n+1)}{2}-1}):=\sum_{i,j}  c_{ij} z_1^{k_{ij}^1}x_1^{m_{ij}^1}(1-x_1^2)^{\frac{n_{ij}^1}{2}}\cdots z_{n-1}^{k_{ij}^{n-1}}x_{n-1}^{m_{ij}^{n-1}}(1-x^2_{n-1})^{\frac{n_{ij}^{n-1}}{2}}\\
	&\cdot(\widetilde{h_{SU(n-1)}})_{ij}(x_{n},\ldots,x_{\frac{n(n-1)}{2}},z_n,\ldots,z_{\frac{n(n+1)}{2}-2}) \;\left(z_{\frac{n(n+1)}{2}-1}\right)^{\frac{l^{n-1}_{ij}}{n-1}}.\nonumber
	\end{align} Here $S^*:=S^1\setminus\{1\}$ to have the function $z^{\frac{1}{n-1}}$ single-valued.
\end{lemma}
The main ingredients of the proof of Lemma \ref{lemma:Main_Lemma} are captured in the following lemma:
\begin{lemma}\label{lemma:combinatorics}
	Let $p,q,k,l\in\N_0$ and $l>0$. Then
	\begin{align*}
	\int_0^{2\pi} e^{i\frac{k}{l}\phi}\;d\phi&=\frac{1}{i}\int_{S^*}z^\frac{k}{l}\frac{dz}{z},
	\end{align*}
	where $S^*:=S^1\setminus\{1\}$ is chosen such that $z^{\frac{k}{l}}$ is analytic on $\C\setminus\R_+$. In addition
	\begin{align*}
	\int_{0}^{\pi/2}\sin^{k+p}(\phi)\cos^{l+q}(\phi)\;d\phi &=\int_0^1 x^{k+p}(1-x^2)^{\frac{l+q-1}{2}}\; dx.
	\end{align*}
\end{lemma}
\begin{proof}
	Both equalities can be found by using a subsitution. The former integral is found by setting $z=e^{i\phi}$ and the latter by $x=\sin(\phi)$.
\end{proof}
\begin{proof}[Proof of Lemma \ref{lemma:Main_Lemma}]
	We use induction on $N$. The $N=2$ case is already proven by M\"uger and Tuset \cite{Mueger}. So assume that the proposition is true for $N-1$. Consider any finite-type function $h$. Then we see that $h^P$ can be expanded by using the multinomial expansion twice:
	\begin{align*}
	\int_G h^P\; dg=\sum_{\sum_{i,j}\beta_{i,j}=P}&\binom{P}{\beta_{1,1},\ldots,\beta_{M,Q}}\int_G \prod_{i,j}\Big(c_{ij}^{\beta_{ij}} e^{i\beta_{ij}k_{ij}^1\phi_1}\sin^{\beta_{ij}m_{ij}^1}(\psi_1)\cos^{\beta_{ij}n_{ij}^1}(\psi_1)\cdots e^{i\beta_{ij}k_{ij}^{N-1}\phi_{N-1}}\\
	&\cdot\sin^{\beta_{ij}m_{ij}^{N-1}}(\psi_{N-1})\cos^{\beta_{ij}n_{ij}^{N-1}}(\psi_{N-1})(h_{SU(N-1)})_{ij}(g_{SU(N-1)})^{\beta_{ij}}e^{i\beta_{ij}l^N_{ij}\omega_{N-1}}\Big)dg.
	\end{align*}
	Filling in the measure given by Equation (\ref{eq:Haar_Measure}) gives
	\begin{align*}
	\int_G h^P\; dg=\sum_{\sum_{i,j}\beta_{i,j}=P}&\binom{P}{\beta_{1,1},\ldots,\beta_{M,Q}}\int_G \prod_{i,j}\left(c_{ij}^{\beta_{ij}} e^{i\beta_{ij}k_{ij}^1\phi_1}\sin^{\beta_{ij}m_{ij}^1}(\psi_1)\cos^{\beta_{ij}n_{ij}^1}(\psi_1)\cdots e^{i\beta_{ij}k_{ij}^{N-1}\phi_{N-1}}\right.\\
	&\left.\cdot\sin^{\beta_{ij}m_{ij}^{N-1}}(\psi_{N-1})\cos^{\beta_{ij}n_{ij}^{N-1}}(\psi_{N-1})(h_{SU(N-1)})_{ij}(g_{SU(N-1)})^{\beta_{ij}}e^{i\beta_{ij}l^N_{ij}\omega_{N-1}}\right)\\
	&\cdot C_N\cos(\psi_{N-1})\sin^{2(N-1)-1}(\psi_{N-1})\left(\prod_{j=1}^{N-2}\sin(\psi_{j})\cos^{2j-1}(\psi_j)d\psi_j d\phi_j\right)\\
	&\cdot d\psi_{N-1} d\phi_{N-1}d\omega_{N-1} dg_{SU(N-1)}\\
	=\sum_{\sum_{i,j}\beta_{i,j}=P}&\left(\binom{P}{\beta_{1,1},\ldots,\beta_{M,Q}}\prod_{i,j}c_{ij}^{\beta_{ij}}\right) \left[\int_{G'}\left(e^{\sum_{i,j}i\beta_{ij}k_{ij}^1\phi_1}\sin^{\sum_{i,j}\beta_{ij}m_{ij}^1}(\psi_1)\cos^{\sum_{i,j}\beta_{ij}n_{ij}^1}(\psi_1)\cdots\right.\right.\\
	&\left.\cdot \, e^{i\sum_{i,j}\beta_{ij}k_{ij}^{N-1}\phi_{N-1}}\sin^{\sum_{i,j}\beta_{ij}m_{ij}^{N-1}}(\psi_{N-1})\cos^{\sum_{i,j}\beta_{ij}n_{ij}^{N-1}}(\psi_{N-1})e^{\sum_{i,j}i\beta_{ij}l^N_{ij}\omega_{N-1}}\right)\\
	&\cdot C_N\cos(\psi_{N-1})\sin^{2(N-1)-1}(\psi_{N-1})\left(\prod_{j=1}^{N-2}\sin(\psi_{j})\cos^{2j-1}(\psi_j)d\psi_j d\phi_j\right)\\
	&\cdot d\psi_{N-1} d\phi_{N-1}d\omega_{N-1}\Bigg]\cdot\left[ \prod_{i,j}\int_{SU(N-1)}(h_{SU(N-1)})_{ij}(g_{SU(N-1)})^{\beta_{ij}} dg_{SU(N-1)}\right],
	\end{align*}
	where we denoted $G'=[0,\pi]\times[0,2\pi]^{N-2}\times[0,\frac{\pi}{2}]^{N-1}\times [0,\frac{2\pi}{N-1}]$ which are the intervals in which $\phi_1,\ldots,\phi_{N-1},\psi_1\ldots,\psi_{N-1}$ and $\omega_{N-1}$ lie, respectively. We note that the integrals over $\omega_{N-1}$ and $\phi_1$ are not over the interval $[0,2\pi]$ yet, hence we will make the substitution $\Omega_{N-1}= (N-1)\omega_{N-1}$ and $\Phi_1=2\phi_1$. Then $d\Omega_{N-1}=(N-1)d\omega_{N-1}$ and $d\Phi_1=2d\phi_1$. This allows us to make use of Lemma \ref{lemma:combinatorics} to rewrite the integral as 

	\begin{align*}
	\int_G h^P\;dg= \sum_{\sum_{i,j}\beta_{i,j}=P}&\frac{C_N}{2(N-1)i^{N}}\left(\binom{P}{\beta_{1,1},\ldots,\beta_{M,Q}}\prod_{i,j}c_{ij}^{\beta_{ij}}\right)\left[\int_{X}\Bigg(z_1^{\sum_{i,j}\beta_{ij}\frac{k_{ij}^1}{2}}x_1^{\sum_{i,j}\beta_{ij}m_{ij}^1}\right.\\
	&\cdot(1-x_1^2)^{\sum_{i,j}\beta_{ij}\frac{n_{ij}^1}{2}}z_{2}^{\sum_{i,j}\beta_{ij}k_{ij}^{2}} x_{2}^{\sum_{i,j}\beta_{ij}m_{ij}^{2}} (1-x_{2}^2)^{\sum_{i,j}\beta_{ij}\frac{n_{ij}^{2}}{2}}\cdots z_{N-1}^{\sum_{i,j}\beta_{ij}k_{ij}^{N-1}}\\
	&\cdot \,x_{N-1}^{\sum_{i,j}\beta_{ij}m_{ij}^{N-1}} (1-x_{N-1}^2)^{\sum_{i,j}\beta_{ij}\frac{n_{ij}^{N-1}}{2}}\left(z_{\frac{N(N+1)}{2}-1}\right)^{\sum_{i,j}\beta_{ij}\frac{l^N_{ij}}{N-1}}\Bigg)\\ 
	&\cdot x_{N-1}^{2N-3}\left(\prod_{j=1}^{N-2}x_j(1-x_j^2)^{j-1}dx_j \frac{dz_j}{z_j}\right)dx_{N-1} \frac{dz_{N-1}}{z_{N-1}}\frac{dz_{\frac{N(N+1)}{2}-1}}{z_{\frac{N(N+1)}{2}-1}}\\
	&\left.\cdot \prod_{i,j}\int_{SU(N-1)}(h_{SU(N-1)})_{ij}(g_{SU(N-1)})^{\beta_{ij}} dg_{SU(N-1)}\right],
	\end{align*}
	where $X=[0,1]^{N-1}\times (S^*)^{N}$. We are now in a position to use the induction hypothesis, which reduces the integral over $h_{SU(N-1)}$ to the following:
	\begin{align*}
	\int_G h^P\;dg= \sum_{\sum_{i,j}\beta_{i,j}=P}&\frac{C_N}{2(N-1)i^{\frac{N(N+1)}{2}-1}}\left(\binom{P}{\beta_{1,1},\ldots,\beta_{M,Q}}\prod_{i,j}c_{ij}^{\beta_{ij}}\right)\left[\int_{X}\Bigg(z_1^{\sum_{i,j}\beta_{ij}\frac{k_{ij}^1}{2}}x_1^{\sum_{i,j}\beta_{ij}m_{ij}^1}\right.\\
	&\cdot(1-x_1^2)^{\sum_{i,j}\beta_{ij}\frac{n_{ij}^1}{2}}z_{2}^{\sum_{i,j}\beta_{ij}k_{ij}^{2}} x_{2}^{\sum_{i,j}\beta_{ij}m_{ij}^{2}} (1-x_{2}^2)^{\sum_{i,j}\beta_{ij}\frac{n_{ij}^{2}}{2}}\cdots z_{N-1}^{\sum_{i,j}\beta_{ij}k_{ij}^{N-1}}\\
	&\cdot x_{N-1}^{\sum_{i,j}\beta_{ij}m_{ij}^{N-1}} (1-x_{N-1}^2)^{\sum_{i,j}\beta_{ij}\frac{n_{ij}^{N-1}}{2}}\left(z_{\frac{N(N+1)}{2}-1}\right)^{\sum_{i,j}\beta_{ij}\frac{l^N_{ij}}{N-1}}\Bigg)\\ 
	&\cdot x_{N-1}^{2N-3}\left(\prod_{j=1}^{N-2}x_j(1-x_j^2)^{j-1}dx_j \frac{dz_j}{z_j}\right)dx_{N-1} \frac{dz_{N-1}}{z_{N-1}}\frac{dz_{\frac{N(N+1)}{2}-1}}{z_{\frac{N(N+1)}{2}-1}}\\
	&\cdot \prod_{i,j}\int_{[0,1]^{\frac{(N-1)(N-2)}{2}}}\int_{(S^*)^{\frac{N(N-1)}{2}-1}}((\widetilde{h_{SU(N-1)}})_{ij}(x_N,\ldots,z_{\frac{N(N+1)}{2}-2}))^{\beta_{ij}}\\
	&\left.\cdot J_{SU(N-1)}(x_N,\ldots,x_{\frac{N(N-1)}{2}})\frac{dz_N}{z_N}\ldots dx_{\frac{N(N-1)}{2}}\right].
	\end{align*}
	Note that we can pull the factors $\beta_{ij}$ back out, which gives
	\begin{align*}
	\int_G h^P\;dg = \sum_{\sum_{i,j}\beta_{i,j}=P}&\frac{1}{2(N-1)i^{\frac{N(N+1)}{2}-1}}\binom{P}{\beta_{1,1},\ldots,\beta_{M,Q}}\int_{[0,1]^{\frac{N(N-1)}{2}}}\int_{S^{\frac{N(N+1)}{2}-1}}\prod_{i,j}\left[c_{ij}z_1^{\frac{k_{ij}^1}{2}}x_1^{m_{ij}^1}\right.\\
	&\cdot (1-x_1^2)^{\frac{n_{ij}^1}{2}}z_{2}^{k_{ij}^{2}}  x_{2}^{m_{ij}^{2}}(1-x_{2}^2)^{\frac{n_{ij}^{2}}{2}}\cdots z_{N-1}^{k_{ij}^{N-1}}x_{N-1}^{m_{ij}^{N-1}} (1-x_{N-1}^2)^{\frac{n_{ij}^{N-1}}{2}}\\ 
	&\left.\cdot\, (\widetilde{h_{SU(N-1)}})_{ij}(x_N,\ldots,z_{\frac{N(N+1)}{2}-2})\left(z_{\frac{N(N+1)}{2}-1}\right)^{\frac{l^N_{ij}}{N-1}}\right]^{\beta_{ij}}\cdot J_{SU(N)}(x_1,\ldots,x_{\frac{N(N-1)}{2}})\\
	&\cdot \frac{dz_1}{z_1}\ldots\frac{dz_{\frac{N(N+1)}{2}-1}}{z_{\frac{N(N+1)}{2}-1}}dx_1\ldots dx_{\frac{N(N-1)}{2}}\\
	=\frac{1}{2(N-1)i^{\frac{N(N+1)}{2}-1}}&\int_{[0,1]^{\frac{N(N-1)}{2}}}\int_{(S^*)^{\frac{N(N+1)}{2}-1}}[\widetilde{h_{SU(N)}}]^P J_{SU(N)}(x_1,\ldots,x_{\frac{N(N-1)}{2}})\frac{dz_1}{z_1}\ldots dx_{\frac{N(N-1)}{2}},
	\end{align*} which is the desired result.
\end{proof}

In other words, we have translated the problem of the non-abelian group $SU(N)$ to the simpler set $[0,1]^{\frac{N(N-1)}{2}}\times (S^*)^{\frac{N(N+1)}{2}-1}$. This is used to translate Mathieu's conjecture to a complex analysis question in the case of $SU(N).$

\begin{definition}
	Let $k,l\in\N$ and $f: [0,1]^{k}\times (S^*)^{l}\rightarrow\C$. We say $f$ is a $SU(N)$-\emph{admissible function} if $f$ can be written as $$f(x_1,\ldots,x_k,z_1,\ldots,z_{l})=\sum_{\vec{m}}c_{\vec{m}}(x)z^{\vec{m}},$$ where $\vec{m}=(m_1,\ldots,m_l)$ is a multi-index where $m_i\in \bigcup_{j=1}^N\frac{1}{j}\mathbb{Z}$, and $c_{\vec{m}}(x)\in \C[x_1,(1-x_1^2)^{1/2},\ldots,x_k,(1-x_k^2)^{1/2}]$ is a complex polynomial in $x_i$ and $\sqrt{1-x_i^2}$. We define the collection of $\vec{m}$ for which $c_{\vec{m}}\neq 0$ \emph{the spectrum of} $f$, and it will be denoted by $\mathrm{Sp}(f)$.
\end{definition} 

It is clear that $\widetilde{h_{SU(N)}}$ is a $SU(N)$-admissible function, so we focus on this class of functions. Motivated by \cite{Mueger}, we make the following conjecture:

\begin{conjecture}\label{con:xz-conjecture}
	Let $f:[0,1]^{\frac{N(N-1)}{2}}\times (S^*)^{\frac{N(N+1)}{2}-1}\rightarrow\C$ be a $SU(N)$-admissible function. If $$\int_{[0,1]^{\frac{N(N-1)}{2}}}\int_{(S^*)^{\frac{N(N+1)}{2}-1}}f^P J_{SU(N)} = 0$$ for all $P\in \N$, then $\vec{0}$ does not lie in the convex hull of $\mathrm{Sp}(f)$.
\end{conjecture}

At first sight, this conjecture may seem to have little to do with Mathieu's conjecture. However

\begin{theorem}\label{thm:Mathieu_proven_assuming_XZ-conjecture}
	Assume Conjecture \ref{con:xz-conjecture} is true. Then Mathieu's conjecture is true for $SU(N)$.
\end{theorem}

\begin{proof}
	Let $f,h$ be finite-type functions of $SU(N)$. Then both are of the form of Equation (\ref{eq:K-finite_function}). Assume $\int_{SU(N)}f^P=0$ for all $P\in \N$. By Lemma \ref{lemma:Main_Lemma}, this is equivalent to $$\int_{[0,1]^{\frac{N(N-1)}{2}}}\int_{(S^*)^{\frac{N(N+1)}{2}-1}}\tilde{f}^P J_{SU(N)}=0,$$ where $\tilde{f}$ is defined as in Equation (\ref{eq:h_tilde}). Applying our assumption gives that $0$ does not lie in the convex hull of $\mathrm{Sp}(\tilde{f})$. Let us write 
	\begin{align}\label{eq:finite_type_function_f_Big_theorem}
	f=f_{SU(N)}(g)&=\sum_{j=1}^M\sum_{i=1}^Q c_{ij} e^{ik_{ij}^1\phi_1}\sin^{m_{ij}^1}(\psi_1)\cos^{n_{ij}^1}(\psi_1)\cdots e^{ik_{ij}^{N-1}\phi_{N-1}}\sin^{m_{ij}^{N-1}}(\psi_{N-1})\cos^{n_{ij}^{N-1}}(\psi_{N-1})\\
	&\qquad\qquad\cdot (f_{SU(N-1)})_{ij}(g_{SU(N-1)}) e^{il^{N-1}_{ij}\omega_{N-1}}\nonumber
	\end{align}
	where the subscript $SU(N)$ indicates it is a finite-type function of $SU(N)$, so that $f_{SU(N-1)}$ is a finite-type function of $SU(N-1)$. Note that by Lemma \ref{lemma:Main_Lemma}, 
	\begin{align}\label{eq:Spectrum_looks_like_this}
	\mathrm{Sp}(\tilde{f})=\bigcup_{i,j}\{(k_{ij}^1,\ldots,k_{ij}^{N(N-1)/2},l_{ij}^1,\ldots,l_{ij}^{N-1})\},
	\end{align} where the constants $k_{ij}^p,l_{ij}^q$ are as in Equation (\ref{eq:finite_type_function_f_Big_theorem}).
	
	We need to prove that $\int_{SU(N)}f^Ph=0$ for $P$ large enough. Assume to the contrary that there exists infinitely many $P$ such that $\int_{SU(N)}f^Ph\neq 0$. The goal of the proof is to show that this gives that $0\in\mathrm{Conv}(\mathrm{Sp}(\tilde{f}))$, taking the identity of Equation (\ref{eq:Spectrum_looks_like_this}) into account.		
	Because of the linearity of the integral, it is enough to show this for $h$ being a monomial. So let us write 
	\begin{align*}
	h=h_{SU(N)}(g)&=e^{iK_1\phi_1}\sin^{R_1}(\psi_1)\cos^{S_1}(\psi_1)\cdots e^{iK_{N-1}\phi_{N-1}}\sin^{R_{N-1}}(\psi_{N-1})\cos^{S_{N-1}}(\psi_{N-1})\\
	&\quad\cdot h_{SU(N-1)}(g_{SU(N-1)}) e^{iL_{N-1}\omega_{N-1}}.		
	\end{align*}
	Note that $h_{SU(N-1)}$ is a monomial finite-type function as well. If $\int_{SU(N)}f^Ph\neq 0$, then there is at least one term over which the integral is non-zero. Going through the same calculations as in previous proof, there is a set if integers $\{\beta_{ij}\}_{i,j}$ such that $\sum_{i,j}\beta_{ij}=P$ and such that
	\begin{align}\label{eq:Haar_invariance_solves_Mathieu}
	\begin{split}0\neq&\int_{G'}e^{i(\sum_{i,j}\beta_{ij}k_{ij}^1+K_1)\phi_1}\sin^{\sum_{i,j}\beta_{ij}m_{ij}^1+R_1}(\psi_1)\cos^{\sum_{i,j}\beta_{ij}n_{ij}^1+S_1}(\psi_1)\cdots e^{i(\sum_{i,j}\beta_{ij}k_{ij}^{N-1}+K_{N-1})\phi_{N-1}}\\
	&\cdot \sin^{\sum_{i,j}\beta_{ij}m_{ij}^{N-1}+R_{N-1}}(\psi_{N-1})\cos^{\sum_{i,j}\beta_{ij}n_{ij}^{N-1}+S_{N-1}}(\psi_{N-1})e^{i(\sum_{i,j}\beta_{ij}l^{N-1}_{ij}+L_{N-1})\omega_{N-1}}\\
	&\cdot C_N\cos(\psi_{N-1})\sin^{2(N-1)-1}(\psi_{N-1})\left(\prod_{j=1}^{N-2}\sin(\psi_{j})\cos^{2j-1}(\psi_j)d\psi_j d\phi_j\right)d\psi_{N-1} d\phi_{N-1}d\omega_{N-1}\\
	&\cdot \left[\prod_{i,j}\int_{SU(N-1)}(f_{SU(N-1)})_{ij}(g_{SU(N-1)})^{\beta_{ij}}h_{SU(N-1)}(g_{SU(N-1)}) dg_{SU(N-1)}\right].\end{split}
	\end{align}
	Our goal is to show that the arguments of all the exponential mappings in Equation (\ref{eq:Haar_invariance_solves_Mathieu}) are zero. To do this, we will make use of properties of the Haar measure. Note that $dg$ is left- and right-invariant, meaning that $\int_{SU(N)}f^{P}(gy)h(gy)dg=\int_{SU(N)}f^P(g)h(g)dg=\int_{SU(N)}f^P(yg)h(yg)dg$ for any $y\in SU(N)$. This must limit the possible parameters. The rest of the proof will therefore consist of choosing convenient matrices $y\in SU(N)$ to get restrictions on these parameters which will prove the proposition.
	
	But before finding these matrices specifically, we give a construction on how to continue restricting the relevant parameters of $f_{SU(N-1)}$ when we have found a construction on restricting the parameters $K_1,\ldots,K_{N-1},L_N$. For note that due to the Euler parametrization, see Lemma \ref{lemma:Euler_Angles}, any $g\in SU(N)$ can be written as 
	$$g=x\begin{pmatrix}
	u&0\\
	0&1
	\end{pmatrix}\xi,$$ where $u\in SU(N-1)$, $\xi=e^{\lambda_{N^2-1}\omega_{N-1}}$ for some $\omega_{N-1}\in\R$ and $x=\prod_{2\leq k\leq N}A(k)(\phi_{k-1},\psi_{k-1})$. In the same way $dg_{SU(N)}=dg_{K}\cdot dg_{SU(N-1)}\cdot d\omega_{N-1}$ for some form $dg_{K}$, as can be seen in Lemma \ref{lemma:Euler_Angles} and Lemma \ref{lemma:Haar_Measure} respectively (for more details on $dg_K$ we refer to our proof of Lemma \ref{lemma:Haar_Measure} and \cite{Helgason_DifGeom_Symm_Spaces}). Specifically, $dg_{SU(N-1)}$ is a Haar measure itself, which means that
	\begin{align*}
	\int_G f^P(g)h(g)dg&=\int_{0}^{\frac{2\pi}{N-1}}\int_{SU(N-1)}\int_{G/K} f^P\left(x\begin{pmatrix}
	u&0\\
	0&1
	\end{pmatrix}\xi\right)h\left(x\begin{pmatrix}
	u&0\\
	0&1
	\end{pmatrix}\xi\right)\, dg_{K}\,dg_{SU(N-1)}\,d\omega_{N-1}\\
	&=\int_{0}^{\frac{2\pi}{N-1}}\int_{SU(N-1)}\int_{G/K} f^P\left(x\begin{pmatrix}
	u'u&0\\
	0&1
	\end{pmatrix}\xi\right)h\left(x\begin{pmatrix}
	u'u&0\\
	0&1
	\end{pmatrix}\xi\right) \, dg_{K}\,dg_{SU(N-1)}\,d\omega_{N-1}\\
	&=\int_{0}^{\frac{2\pi}{N-1}}\int_{SU(N-1)}\int_{G/K} f^P\left(x\begin{pmatrix}
	uu'&0\\
	0&1
	\end{pmatrix}\xi\right)h\left(x\begin{pmatrix}
	uu'&0\\
	0&1
	\end{pmatrix}\xi\right) \, dg_{K}\,dg_{SU(N-1)}\,d\omega_{N-1}
	\end{align*}
	for any $u'\in SU(N-1)$ by using properties of the Haar measure $dg_{SU(N-1)}$. Here we denoted $G/K$ as the following space
	\begin{align*}
	G/K=\left\{\left.\prod_{2\leq k\leq N}A(k)(\phi_{k-1},\psi_{k-1})\,\right|\,\phi_k\in[0,2\pi],\psi_k\in\left[0,\frac{\pi}{2}\right]\right\}.
	\end{align*}

	Looking at which parameters change in Equation (\ref{eq:Haar_invariance_solves_Mathieu}) when changing $u$ to $u'u$ or $uu'$, we see that the following equation must hold 
	\begin{align*}
	\begin{split}\int_{SU(N-1)}f_{SU(N-1)}&(g_{SU(N-1)})_{ij}^{\beta_{ij}}h_{SU(N-1)}(g_{SU(N-1)})\,dg_{SU(N-1)}\\
	&=\int_{SU(N-1)}f_{SU(N-1)}(u'g_{SU(N-1)})_{ij}^{\beta_{ij}}h_{SU(N-1)}(u'g_{SU(N-1)})\,dg_{SU(N-1)}\\
	&=\int_{SU(N-1)}f_{SU(N-1)}(g_{SU(N-1)}u')_{ij}^{\beta_{ij}}h_{SU(N-1)}(g_{SU(N-1)}u')\,dg_{SU(N-1)}\end{split}
	\end{align*}
	for any $u'\in SU(N-1)$. This shows that any construction on $SU(N)$ to restrict the parameters $K_1,\ldots,K_{N-1},L_N$ can also be applied to $SU(N-1)$ and the parameters $K_{N},\ldots, K_{\frac{N(N-1)}{2}},L_{1},\ldots,L_{N-2}$ from those finite-type functions, yielding the same result. It is therefore enough to know what the restrictions of $K_1,\ldots,K_{N-1},L_{N-1}$ are.
	
	Now let us define
	\begin{align}
	D_{k,n}(z):=\mathrm{diag}(e^{iz},\ldots,e^{iz},e^{-i(n-1)z},1,\ldots,1),
	\end{align} which is a $k\times k$ matrix, where $2\leq n\leq k$ and $z\in\R$. Here the diagonal has $n-1$ times $e^{iz}$, and $k-n$ ones. Then $D_{k,n}(z)\in SU(k)$ for all $n$ and $z$. Recall that by the properties of the Haar measure, the mapping $g\mapsto D_{N,2}(z)g$ is invariant. That is to say, the map $L_{D_{N,2}(z)}:G\rightarrow G$, given by $L_{D_{N,2}(z)}g=D_{N,2}(z)\,g$, is invariant, i.e. $$\int_G f^P(L_{D_{N,2}(z)}(g)) h(L_{D_{N,2}(z)}(g))dg = \int_G f^P(D_{N,2}(z)\,g)h(D_{N,2}(z)\,g)dg = \int_G f^P(g)h(g)dg.$$
	So note if $g\in SU(N)$ we have
	\begin{align*}
	D_{N,2}(z)g&=\begin{pmatrix}
	e^{iz}&&\\
	&e^{-iz}&\\
	&&\mathbf{1}_{N-2}
	\end{pmatrix} e^{\lambda_3\phi_1}e^{\lambda_2\psi_1}\ldots e^{\lambda_3\phi_{N-1}}e^{\lambda_{(N-1)^2+1}\psi_{N-1}}\begin{pmatrix}
	F_{N-1}(\phi_N,\ldots\omega_{N-2})&0\\
	0&1
	\end{pmatrix}e^{\lambda_{N^2-1}\omega_{N-1}}\\
	&=e^{\lambda_3(\phi_1+z)}e^{\lambda_2\psi_1}\ldots e^{\lambda_3\phi_{N-1}}e^{\lambda_{(N-1)^2+1}\psi_{N-1}}\begin{pmatrix}
	F_{N-1}(\phi_N,\ldots\omega_{N-2})&0\\
	0&1
	\end{pmatrix}e^{\lambda_{N^2-1}\omega_{N-1}}.
	\end{align*}
	This means, by bijectivity of the Euler parametrization that
	$$D_{N,2}(z)g=F_N(\phi_1+z,\phi_2,\ldots,\phi_{\frac{N(N-1)}{2}},\psi_1,\ldots,\psi_{\frac{N(N-1)}{2}},\omega_1,\ldots,\omega_{N-1}),$$	
	which shows that the mapping $g\mapsto D_{N,2}(z)g$ is equivalent to sending $\phi_1\mapsto \phi_1+z$. Since $g\mapsto D_{N,2}(z)g$ is invariant under the Haar measure, this means that sending $\phi_1\mapsto \phi_1+z$ for any $z\in\R$ should be invariant as well in Equation (\ref{eq:Haar_invariance_solves_Mathieu}), which is to say the integral does not change if we replace $\phi$ with $\phi+z$. This can only be the case in Equation (\ref{eq:Haar_invariance_solves_Mathieu}) if 
	\begin{align*}
	\sum_{i,j}\beta_{ij}k_{ij}^1+K_1=0.
	\end{align*} 
	In the same way we see that the mapping $g\mapsto gD_{N,N}(z)$ is equivalent to sending $\omega_{N-1}\mapsto \omega_{N-1}+z$ which should be invariant. This can only lead to the same equation in Equation (\ref{eq:Haar_invariance_solves_Mathieu}) if 
	\begin{align}\label{eq:Parameters_SU(N)_L_N_conv_hull}
	\sum_{i,j}\beta_{ij}l_{ij}^{N-1}+L_{N-1}=0.
	\end{align}
	As stated before, this construction can be applied on $f_{SU(N-1)}^Ph_{SU(N-1)}$ as well with the matrix $D_{N-1,n}(z)$. Setting $n=2$ and $n=N-1$ give the equations 
	\begin{align}\label{eq:Parameters_SU(N-1)_conv_hull}
	\sum_{i,j}\beta_{ij}k_{ij}^N+K_N=0,\qquad\sum_{ij}\beta_{ij}l_{ij}^{N-1}+L_{N-1}=0.
	\end{align}
	Next, we note that we can write the measure as a product of the Haar measure on $SU(2)$ and other measures in the following way:
	\begin{align*}
	dg_{SU(N)}=T(\psi_2,\ldots,\psi_{N-1})  dg_{SU(2)}d\psi_2 d\psi_3d\phi_3\ldots d\psi_{N-1}d\phi_{N-1}d\omega_{N-1}dg_{SU(N-1)},
	\end{align*}
	where $T(\psi_2,\ldots,\psi_{N-1})=2\pi^2C_N \cos(\psi_{N-1})\sin^{2(N-1)-1}(\psi_{N-1})\left[\prod_{j=2}^{N-2}\sin(\psi_{j})\cos^{2j-1}(\psi_j)\right]$. The Euler para\-metrization also allows us to write $$g_{SU(N)}=F_2(\phi_1,\psi_1,\phi_2)e^{\lambda_{5}\psi_{2}}\left(\prod_{3\leq k\leq N}A(k)(\phi_{k-1},\psi_{k-1})\right)\begin{pmatrix}
	F_{N-1}(\phi_N,\ldots,\omega_{N-2})&0\\
	0&1
	\end{pmatrix}e^{\lambda_{N^2-1}\omega_{N-1}}.$$ This way we have found a subset isomorphic to $SU(2)$ over which we integrate with the Haar measure associated with $SU(2)$ group. This means that we can use the features of the Haar measure on $SU(2)$, which is to say that for any function $H:SU(N)\rightarrow\C$ we have the equation
	\begin{align*}
	&\int_{SU(N)} H\left(e^{\lambda_3\phi_1}e^{\lambda_2\psi_1}e^{\lambda_3\phi_2}e^{\lambda_5\psi_2}\ldots e^{\lambda_{(N-1)^2+1}\psi_{N-1}}\begin{pmatrix}
	F_{N-1}(\phi_N,\ldots\omega_{N-2})&0\\
	0&1
	\end{pmatrix}e^{\lambda_{N^2-1}\omega_{N-1}}\right) \,dg\\
	&= \int_{SU(N)} H\left(\begin{pmatrix}
	A&0\\
	0&\mathbf{1}_{N-2}
	\end{pmatrix}e^{\lambda_3\phi_1}e^{\lambda_2\psi_1}e^{\lambda_3\phi_2}e^{\lambda_5\psi_2}\ldots e^{\lambda_{(N-1)^2+1}\psi_{N-1}}\begin{pmatrix}
	F_{N-1}(\phi_N,\ldots\omega_{N-2})&0\\
	0&1
	\end{pmatrix}e^{\lambda_{N^2-1}\omega_{N-1}}\right) \,dg\\
	&= \int_{SU(N)} H\left(e^{\lambda_3\phi_1}e^{\lambda_2\psi_1}e^{\lambda_3\phi_2}\begin{pmatrix}
	A&0\\
	0&\mathbf{1}_{N-2}
	\end{pmatrix}e^{\lambda_5\psi_2}\ldots e^{\lambda_{(N-1)^2+1}\psi_{N-1}}\begin{pmatrix}
	F_{N-1}(\phi_N,\ldots\omega_{N-2})&0\\
	0&1
	\end{pmatrix}e^{\lambda_{N^2-1}\omega_{N-1}}\right) \,dg
	\end{align*}
	for any $A\in SU(2)$. In specific, taking $A=\begin{pmatrix}
	e^{iz}&0\\
	0&e^{-iz}
	\end{pmatrix}$ for any $z\in \R$ gives that the following map is also invariant:
	\begin{align*}
	&e^{\lambda_3\phi_1}e^{\lambda_2\psi_1}e^{\lambda_3\phi_2}e^{\lambda_5\psi_2}\ldots e^{\lambda_3\phi_{N-1}}e^{\lambda_{(N-1)^2+1}\psi_{N-1}}\begin{pmatrix}
	F_{N-1}(\phi_N,\ldots\omega_{N-2})&0\\
	0&1
	\end{pmatrix}e^{\lambda_{N^2-1}\omega_{N-1}}\\
	&\mapsto e^{\lambda_3\phi_1}e^{\lambda_2\psi_1}e^{\lambda_3\phi_2}\begin{pmatrix}
	e^{iz}&&\\
	&e^{-iz}&\\
	&&\mathbf{1}_{N-2}
	\end{pmatrix}e^{\lambda_5\psi_2}\ldots e^{\lambda_3\phi_{N-1}}e^{\lambda_{(N-1)^2+1}\psi_{N-1}}\begin{pmatrix}
	F_{N-1}(\phi_N,\ldots\omega_{N-2})&0\\
	0&1
	\end{pmatrix}e^{\lambda_{N^2-1}\omega_{N-1}}\\
	&=e^{\lambda_3\phi_1}e^{\lambda_2\psi_1}e^{\lambda_3(\phi_2+z)}e^{\lambda_5\psi_2}\ldots e^{\lambda_3\phi_{N-1}}e^{\lambda_{(N-1)^2+1}\psi_{N-1}}\begin{pmatrix}
	F_{N-1}(\phi_N,\ldots\omega_{N-2})&0\\
	0&1
	\end{pmatrix}e^{\lambda_{N^2-1}\omega_{N-1}}.
	\end{align*}
	Hence sending $\phi_2\mapsto \phi_2+z$ is also an invariance for all $z\in\R$. Looking back at Equation (\ref{eq:Haar_invariance_solves_Mathieu}) this can only be the case if
	\begin{align*}
	\sum_{i,j}\beta_{ij}k_{ij}^2+K_2=0.
	\end{align*}
	Finally, to get an equation for $K_3,\ldots,K_{N-1}$, we see that
	\begin{align*}
	&D_{N,n}(z)\begin{pmatrix}
	\cos(\psi_{n-1})&0&\ldots&0&\sin(\psi_{n-1})&&&\\
	0&1&&&0&&&\\
	\vdots&&\ddots&&\vdots&&&\\
	0&&&1&0&&&\\
	-\sin(\psi_{n-1})&0&\ldots&0&\cos(\psi_{n-1})&&&\\
	&&&&&1&&\\
	&&&&&&\ddots&\\
	&&&&&&&1
	\end{pmatrix}=\\
	&\begin{pmatrix}
	e^{inz}&&\\
	&e^{-inz}&\\
	&&\mathbf{1}_{N-2}
	\end{pmatrix}		\begin{pmatrix}
	\cos(\psi_{n-1})&0&\ldots&0&\sin(\psi_{n-1})&&&\\
	0&1&&&0&&&\\
	\vdots&&\ddots&&\vdots&&&\\
	0&&&1&0&&&\\
	-\sin(\psi_{n-1})&0&\ldots&0&\cos(\psi_{n-1})&&&\\
	&&&&&1&&\\
	&&&&&&\ddots&\\
	&&&&&&&1
	\end{pmatrix}\begin{pmatrix}
	e^{-inz}&&\\
	&e^{inz}&\\
	&&\mathbf{1}_{N-2}
	\end{pmatrix}D_{N,n}(z),
	\end{align*}
	where the block matrix is an $n\times n$ matrix. In addition, if the block matrix is an $m\times m$ matrix with $m>n$, then we see that 
	\begin{align*}
	&D_{N,n}(z)\begin{pmatrix}
	\cos(\psi_{m-1})&0&\ldots&0&\sin(\psi_{m-1})&&&\\
	0&1&&&0&&&\\
	\vdots&&\ddots&&\vdots&&&\\
	0&&&1&0&&&\\
	\sin(\psi_{m-1})&0&\ldots&0&\cos(\psi_{m-1})&&&\\
	&&&&&1&&\\
	&&&&&&\ddots&\\
	&&&&&&&1
	\end{pmatrix}=\\
	&\begin{pmatrix}
	e^{iz}&&\\
	&e^{-iz}&\\
	&&\mathbf{1}_{N-2}
	\end{pmatrix}		\begin{pmatrix}
	\cos(\psi_{m-1})&0&\ldots&0&\sin(\psi_{m-1})&&&\\
	0&1&&&0&&&\\
	\vdots&&\ddots&&\vdots&&&\\
	0&&&1&0&&&\\
	\sin(\psi_{m-1})&0&\ldots&0&\cos(\psi_{m-1})&&&\\
	&&&&&1&&\\
	&&&&&&\ddots&\\
	&&&&&&&1
	\end{pmatrix}\begin{pmatrix}
	e^{-iz}&&\\
	&e^{iz}&\\
	&&\mathbf{1}_{N-2}
	\end{pmatrix}D_{N,n}(z).
	\end{align*}
	Given that $D_{N,n}(z)$ is identical to $e^{iz}\mathbf{1}_n$ in the upper left $(n-1)\times (n-1)$ corner, it commutes with $e^{\lambda_3\phi_k}$ and $e^{\lambda_{(l-1)^2+1}\lambda_{l-1}}$ for all $k$ and $l\leq n-1$. This way, we see that
	\begin{align*}
	D_{N,N}(z)g&= D_{N,N}(z)e^{\lambda_3\phi_1}e^{\lambda_2\psi_1}\ldots e^{\lambda_3\phi_{N-1}}e^{\lambda_{(N-1)^2+1}\psi_{N-1}}\begin{pmatrix}
	F_{N-1}(\phi_N,\ldots\omega_{N-2})&0\\
	0&1
	\end{pmatrix}e^{\lambda_{N^2-1}\omega_{N-1}}\\
	&=e^{\lambda_3\phi_1}e^{\lambda_2\psi_1}\ldots e^{\lambda_3\phi_{N-1}}D_{N,N}(z)e^{\lambda_{(N-1)^2+1}\psi_{N-1}}\begin{pmatrix}
	F_{N-1}(\phi_N,\ldots\omega_{N-2})&0\\
	0&1
	\end{pmatrix}e^{\lambda_{N^2-1}\omega_{N-1}}\\
	&=e^{\lambda_3\phi_1}e^{\lambda_2\psi_1}\ldots e^{\lambda_3(\phi_{N-1}+zn)}e^{\lambda_{(N-1)^2+1}\psi_{N-1}}e^{-nz\lambda_3}D_{N,N}(z)\begin{pmatrix}
	F_{N-1}(\phi_N,\ldots\omega_{N-2})&0\\
	0&1
	\end{pmatrix}e^{\lambda_{N^2-1}\omega_{N-1}}\\
	&=e^{\lambda_3\phi_1}e^{\lambda_2\psi_1}\ldots e^{\lambda_3(\phi_{N-1}+zn)}e^{\lambda_{(N-1)^2+1}\psi_{N-1}}e^{-nz\lambda_3}\begin{pmatrix}
	F_{N-1}(\phi_N,\ldots\omega_{N-2})&0\\
	0&1
	\end{pmatrix}e^{\lambda_{N^2-1}(\omega_{N-1}+z)}.
	\end{align*}
	We note that $e^{-nz\lambda_3}\in \begin{pmatrix}
	SU(N-1)&0\\
	0&1
	\end{pmatrix}$, and by previous arguments the mapping $F_{N-1}(\phi_N,\ldots,\omega_{N-2})\mapsto e^{-nz\lambda_3}F_{N-1}(\phi_N,\ldots,\omega_{N-2})$ is equivalent to $\phi_N\mapsto \phi_N-nz$. So in total we see that the mapping $g\mapsto D_N(z)g$ is equivalent to the mapping $(\phi_{N-1},\phi_N,\omega_{N-1})\mapsto (\phi_{N-1}+zn,\phi_N-zn,\omega_{N-1}+z).$ Looking back at Equation (\ref{eq:Haar_invariance_solves_Mathieu}) this can only be invariant if
	\begin{align*}
	\left(\sum_{i,j}\beta_{ij}k_{ij}^{N-1}+K_{N-1}\right)+\left(\sum_{ij}\beta_{ij}k_{ij}^N+K_N\right)+\left(\sum_{i,j}\beta_{ij}l_{ij}^{N-1}+L_{N-1}\right)=0.
	\end{align*}
	Combining this with Equations $(\ref{eq:Parameters_SU(N)_L_N_conv_hull})$ and $(\ref{eq:Parameters_SU(N-1)_conv_hull})$, we immediately find
	\begin{align}\label{eq:Parameters_SU(N)_K_N-1_conv_hull}
	\sum_{i,j}\beta_{ij}k_{ij}^{N-1}+K_{N-1}=0.
	\end{align}
	Next, we consider $g\mapsto D_{N,N-1}(z)g$. We then see
	\begin{align*}
	D_{N,N-1}(z)g&= D_{N,N-1}(z)e^{\lambda_3\phi_1}e^{\lambda_2\psi_1}\ldots e^{\lambda_3\phi_{N-1}}e^{\lambda_{(N-1)^2+1}\psi_{N-1}}\begin{pmatrix}
	F_{N-1}(\phi_N,\ldots\omega_{N-2})&0\\
	0&1
	\end{pmatrix}e^{\lambda_{N^2-1}\omega_{N-1}}\\
	&= e^{\lambda_3\phi_1}e^{\lambda_2\psi_1}\ldots e^{\lambda_3\phi_{N-2}}D_{N,N-1}(z)e^{\lambda_{(N-2)^2+1}\psi_{N-2}}e^{\lambda_3\phi_{N-1}}e^{\lambda_{(N-1)^2+1}\psi_{N-1}}\\
	&\quad\cdot\begin{pmatrix}
	F_{N-1}(\phi_N,\ldots\omega_{N-2})&0\\
	0&1
	\end{pmatrix}e^{\lambda_{N^2-1}\omega_{N-1}}\\
	&= e^{\lambda_3\phi_1}e^{\lambda_2\psi_1}\ldots e^{\lambda_3(\phi_{N-2}+(N-1)z)}e^{\lambda_{(N-2)^2+1}\psi_{N-2}}e^{\lambda_{3}(\phi_{N-1}-(N-1)z)}D_{N,N-1}(z)e^{\lambda_{(N-1)^2+1}\psi_{N-1}}\\
	&\quad\cdot \begin{pmatrix}
	F_{N-1}(\phi_N,\ldots\omega_{N-2})&0\\
	0&1
	\end{pmatrix}e^{\lambda_{N^2-1}\omega_{N-1}}\\
	&=e^{\lambda_3\phi_1}e^{\lambda_{2}\psi_1}\ldots e^{\lambda_3\phi_{N-3}}e^{\lambda_{(N-3)^2+1}\psi_{N-3}}e^{\lambda_{3}(\phi_{N-2}+(N-1)z)}e^{\lambda_{(N-2)^2+1}\psi_{N-2}}e^{\lambda_3(\phi_{N-1}-(N-2)z)}\\
	&\quad\cdot e^{\lambda_{(N-1)^2+1}\psi_{N-1}}\begin{pmatrix}
	D_{N-1,2}(-z)D_{N-1,N-1}(z)F_{N-1}(\phi_N,\ldots\omega_{N-2})&0\\
	0&1
	\end{pmatrix}e^{\lambda_{N^2-1}\omega_{N-1}}.
	\end{align*}
	Therefore we see that the mapping $g\mapsto D_{N,N-1}(z)g$ is equivalent to sending $$(\phi_{N-2},\phi_{N-1},F_{N-1}(\phi_{N},\ldots,\omega_{N-2}))\mapsto (\phi_{N-2}+(N-1)z, \phi_{N-1}-(N-2)z,D_{N-1,2}(-z)D_{N-1,N-1}(z)F_{SU(N-1)}).$$
	Note that this transformation contains the mapping $g'\mapsto D_{N-1,2}(-z)D_{N-1,N-1}(z)g'$ where $g'\in SU(N-1)$. By previous arguments how to apply our strategy to $SU(N-1)$, we see that this transformation has already been considered and as a result we could put some parameters equal to 0. So we can ignore transformation $g'\mapsto D_{N-1,N-1}(-z)D_{N-1,2}(z)g'$. In other words, the mapping $g\mapsto D_{N,N-1}(z)g$ is equivalent to the mappings $$(\phi_{N-2},\phi_{N-1})\mapsto (\phi_{N-2}+(N-1)z, \phi_{N-1}-(N-2)z)$$ to see what kind of restrictions we can put on the parameters. Looking at Equation (\ref{eq:Haar_invariance_solves_Mathieu}), this is can only happen if
	\begin{align*}
	(N-1)\left(\sum_{i,j}\beta_{ij}k_{ij}^{N-2}+K_{N-2}\right)+(N-2)\left(\sum_{i,j}\beta_{ij}k_{ij}^{N-1}+K_{N-1}\right)=0.
	\end{align*}
	Combining this with Equation (\ref{eq:Parameters_SU(N)_K_N-1_conv_hull}) we see that $$\sum_{i,j}\beta_{ij}k_{ij}^{N-2}+K_{N-2}=0.$$ This process can be repeated for $g\mapsto D_{N-m}g$ with $m=2,\ldots,N-3$ increasing, and in the end we find that
	\begin{align*}
	\sum_{i,j}\beta_{ij}k_{ij}^{N-m-1}+K_{N-m-1}=0\qquad\forall m\in\{2,\ldots,N-3\}.
	\end{align*}
	Combining everything, we thus see that
	\begin{align*}
	\sum_{i,j}\beta_{ij}k_{ij}^{m}+K_{m}=0\qquad\forall m\in\{1,\ldots,N-1\}.
	\end{align*}
	As discussed, this procedure can be continued on the remaining parameters of $f_{SU(N-1)}$, which gives in the end 
	\begin{align*}
	\sum_{i,j}\beta_{ij}k_{ij}^m+K_m=0,\qquad \sum_{i,j}\beta_{ij}l_{ij}^{n}+L_n=0,\qquad\forall m\in\left\{1,2,\ldots,\frac{N(N-1)}{2}\right\},\, n\in\{1\ldots,N\}.
	\end{align*}	
	Therefore, we see that $$(K_1,K_2,\ldots,K_{N(N-1)/2},L_1,\ldots,L_N)=-\sum_{ij}(k_{ij}^1,k^2_{ij},\ldots,k^{N(N-1)/2}_{ij},l_{ij}^1,\ldots,l_{ij}^{N-1})\beta_{ij}.$$
	Thus we see that 
	$$\frac{1}{P}(K_1,K_2,\ldots,K_{N(N-1)/2},L_1,\ldots,L_N)=-\sum_{ij}(k_{ij}^1,k^2_{ij},\ldots,k^{N(N-1)/2}_{ij},l_{ij}^1,\ldots,l_{ij}^{N-1})\frac{\beta_{ij}}{P}.$$ As $P\rightarrow\infty$, the left-hand side tends to 0. Hence there exists a limiting sequence in the convex hull of $\mathrm{Sp}(\tilde{f})$ that converges to 0. Since the convex hull is a closed set, it shows that $0\in \mathrm{Conv}(\mathrm{Sp}(\tilde{f}))$ which is a contradiction with our assumption. Therefore it cannot be that $\int_{SU(N)}f^Ph\neq0$ for infinitely many $P$, and thus $\int_{SU(N)}f^Ph=0$ for large $P$ as was to be proven.	
\end{proof}

\begin{remark}
	One could wonder whether the inverse implication of Theorem \ref{thm:Mathieu_proven_assuming_XZ-conjecture} is true as well. It is still an open question whether this is the case.
\end{remark}

\section{The case of $SO(N)$}\label{sec:SO(N)}

Having considered the group $SU(N)$, we employ a similar strategy for $SO(N)$. We will leave the details to the reader because of the similarity. Since $SO(N)\subset SU(N)$ is a closed subgroup, we could describe $\mathfrak{so}(N)$ as a subalgebra of $\mathfrak{su}(N)$. We note that the set $$\{\lambda_{j^2-1+2k}|j=1,\ldots N-1, k=1,\ldots,2j\}$$ is a basis of $\mathfrak{so}(N)$. For the rest of this paper we will use this basis to describe $\mathfrak{so}(N)$, and with it we can describe the generalized Euler angles for $SO(n)$:
\begin{restatable}[Generalized Euler Angles on $SO(N)$]{lemma}{EulerAngles_SO(N)}\label{lemma:Euler_Angles_SO(N)}
	Let $N\geq 2$. Define inductively the mapping $\Phi_N:([0,2\pi]\times[0,\pi]^{N-2})\times([0,2\pi]\times[0,\pi]^{N-3})\times\cdots\times([0,2\pi]\times [0,\pi])\times[0,2\pi]\rightarrow SO(N)$ by $\Phi_1\equiv 1$ and
	\begin{equation}\label{eq:Euler_parametrization_SO(N)}
	\Phi_N(\phi_1,\ldots\phi_{\frac{N(N-1)}{2}}):=
	\left(\prod_{1\leq k\leq N-1}e^{\phi_k \lambda_{(k+1)^2-2}}\right)\cdot\begin{pmatrix}
	\Phi_{N-1}(\phi_{N},\ldots,\phi_{\frac{N(N-1)}{2}})&0\\
	0&1
	\end{pmatrix},	
	\end{equation}
	where we denoted the product as $$\prod_{1\leq k\leq N-1}e^{\phi_k \lambda_{(k+1)^2-2}}:=e^{\phi_1\lambda_2}\cdot \cdots\cdot e^{\phi_{N-1}\lambda_{N^2-2}}.$$ This mapping is surjective. Moreover it is a diffeomorphism on the interior of the hypercube.
\end{restatable}

\begin{restatable}{lemma}{HaarMeasure_SO(N)}\label{lemma:Haar_Measure_SO(N)}
	Let $N\geq 2$ and $\Phi_N$ be the Euler parametrization of $SO(N)$ as in Lemma \ref{lemma:Euler_Angles_SO(N)}. The normalised Haar measure $dg_{SO(N)}$ with this parametrization is then given inductively by $$dg_{SO(2)}=C_2\,d\phi_1$$ and
	\begin{align}\label{eq:Haar_Measure_SO(N)}
	\begin{split}dg_{SO(N)}=C_N\prod_{j=1}^{N-1}\sin^{j-1}(\phi_j)d\phi_1\ldots d\phi_{N-1}dg_{SO(N-1)},\end{split}
	\end{align}
	where we denoted $\sin^0(\phi_k)=1$ and $C_n:= \frac{1}{2\pi^{\frac{n}{2}}}\Gamma(\frac{n}{2})$ for all $n=2,\ldots, N$, where $\Gamma$ is the Euler Gamma function.
\end{restatable}
As stated before, the proofs proceed analogously to the proofs given in the appendix, hence will be omitted. Next we will describe the finite-type functions just as in previous section. As in the case of $SU(N)$, the finite-type functions $f$ are sums of products of matrix coefficients since the irreducible representations of $SO(N)$ are polynomials in the matrix coefficients. By the parametrization given in Equation (\ref{eq:Euler_parametrization_SO(N)}), we see that these products consist of (powers of) $\sin(\phi_k), \cos(\phi_k)$, $\exp(i\phi_{\frac{(N-1)N}{2}+1-\frac{(N-k)(N-k+1)}{2}})$ and finite-type functions on $SO(N-1)$ where $k=1,\ldots,N-1$. Again, using $\cos^2(\psi_k)+\sin^2(\psi_k)=1$, any finite-type function $h$ can be written as
\begin{align}\begin{split}\label{eq:K-finite_function_SO(N)}
h_{SO(N)}(g)= \sum_{j=1}^M\sum_{i=1}^Q c_{ij}& e^{ik_{ij}^1\phi_1}\cos^{m_{ij}^1}(\phi_2)\sin^{n_{ij}^1}(\phi_2)\cdots \cos^{m_{ij}^{N-1}}(\phi_{N-1})\sin^{n_{ij}^{N-1}}(\phi_{N-1})\\
&\cdot (h_{SO(N-1)})_{ij}(g_{SO(N-1)}),\end{split}
\end{align}
where $g_{SO(N-1)}:=\Phi_{N-1}(\phi_{N},\ldots,\phi_{\frac{N(N-1)}{2}})$ is the $SO(N-1)$ component of $g=\Phi_N(\phi_1\ldots\phi_{N(N-1)/2})$ as in Lemma \ref{lemma:Euler_Angles_SO(N)}, and $(h_{SO(N-1)})_{ij}$ is a family of finite-type functions on $SO(N-1)$. Again, $k_{ij}^p,l_{ij}^p\in\mathbb{Z}$, $m_{ij}^p\in\N$ and $n_{ij}^p\in\{0,1\}$. We sum over both $i$ and $j$. The sum over $i$ is to ensure we have all possible combinations of different terms, while the sum over $j$ allows for different powers of each term. 

The $SO(N)$ finite-type functions differ slightly from the $SU(N)$ finite-type functions in Equation (\ref{eq:K-finite_function}). We chose to have (possible) higher powers of $\cos(\phi_j)$ instead of $\sin(\phi_j)$ because the Haar measure only contains powers of $\sin(\phi_j)$. In addition, there are fewer parameters going over $[0,2\pi]$, hence we can only write $\phi_1,\phi_N,\ldots,\phi_{N(N-1)/2}$ as an exponential.

In the same way as with $SU(N)$, we can translate the problem then back to analysis of functions on $\R^n\times \C^m$ in the following way:

\begin{lemma}\label{lemma:Main_Lemma_SO(N)}
	Let $h_{SO(N)}$ be a finite-type function on $SO(N)$, as in Equation (\ref{eq:K-finite_function_SO(N)}) and $N\geq 2$. Then for any $P\in \N$ we have
	\begin{align}
	\begin{split}\int_{SO(N)}h(g)^P\,dg=\frac{1}{i^{N-1}}&\int_{[-1,1]^{\frac{(N-1)(N-2)}{2}}}\int_{(S^1)^{N-1}}\left[\widetilde{h_{SO(N)}}(x_1,\ldots,z_{n-1})\right]^P\\
	&\cdot J_{SO(N)}(x_1,\ldots,x_{\frac{(N-1)(N-2)}{2}})\frac{dz_1}{z_1}\ldots\frac{dz_{N-1}}{z_{N-1}}dx_1\ldots dx_{\frac{(N-1)(N-2)}{2}}.\end{split}
	\end{align}
	Here $J_{SO(N)}$ is defined recursively by $J_{SO(2)}\equiv 1$ and, for $3\leq n\leq N$, by $$J_{SO(n)}(x_1,\ldots,x_{\frac{(n-1)(n-2)}{2}}):=C_n\left(\prod_{j=1}^{n-1}(1-x_j^2)^{\frac{j-1}{2}}\right) J_{SO(n-1)}\left(x_{n},\ldots,x_{\frac{(n-1)(n-2)}{2}}\right),$$ where $C_n$ is defined as in Lemma \ref{lemma:Haar_Measure_SO(N)} and $\widetilde{h_{SO(N)}}$ is defined recursively by $\widetilde{h_{SO(1)}}=1$ and by 
	\begin{align}
	\widetilde{h_{SO(n)}}(x_1,\ldots x_{\frac{(n-1)(n-2)}{2}},z_1,\ldots z_{n-1}):=\sum_{i,j} & c_{ij} z_1^{k_{ij}^1}x_1^{m_{ij}^1}(1-x_1^2)^{\frac{n_{ij}^1}{2}}\cdots x_{n-1}^{m_{ij}^{n-1}}(1-x^2_{n-1})^{\frac{n_{ij}^{n-1}}{2}}\\
	&\cdot(\widetilde{h_{SO(n-1)}})_{ij}(x_n,\ldots,x_{\frac{(n-1)(n-2)}{2}},z_{n-2},\ldots,z_{n-1})\nonumber.
	\end{align}
\end{lemma}

Note that Lemma \ref{lemma:Main_Lemma_SO(N)} is similar to Lemma \ref{lemma:Main_Lemma}, the difference here being that the $x_j$ variables go over the interval $[-1,1]$ instead of $[0,1]$ which is due to the original intervals being $[0,\pi]$ instead of $[0,\pi/2]$ in the $SU(N)$ case. The proof goes identical to the proof of Lemma \ref{lemma:Main_Lemma}.

Similar to previous section, we describe a conjecture that will deal with Mathieu's conjecture for $SO(N)$:

\begin{definition}
	Let $k,l\in \N$ and $f: [-1,1]^{k}\times (S^1)^{l}\rightarrow\C$. We say $f$ is a \emph{$SO(N)$-admissible function} if $f$ can be written as $$f(x_1,\ldots,x_k,z_1,\ldots,z_{l})=\sum_{\vec{m}}c_{\vec{m}}(x)z^{\vec{m}},$$ where $\vec{m}=(m_1,\ldots,m_l)\in \mathbb{Z}^l$ is a multi-index, and $c_{\vec{m}}(x)\in \C[x_1,(1-x_1^2)^{1/2},\ldots,x_k,(1-x_k^2)^{1/2}]$ is a complex polynomial in $x_i$ and $\sqrt{1-x_i^2}$. We will call the collection of $\vec{m}$ for which $c_{\vec{m}}\neq 0$ the \emph{spectrum of $f$} and will be denoted by $\mathrm{Sp}(f)$.
\end{definition} 
\begin{conjecture}\label{con:xz-conjecture_SO(N)}
	Let $f:[-1,1]^{\frac{(N-1)(N-2)}{2}}\times S^{N-1}\rightarrow\C$ be a $SO(N)$-admissible function. If $$\int_{[-1,1]^{\frac{(N-1)(N-2)}{2}}}\int_{S^{N-1}}f^P J_{SO(N)} = 0$$ for all $P\in \N$, then $\vec{0}$ does not lie in the convex hull of $\mathrm{Sp}(f)$.
\end{conjecture}

\begin{proposition}
	Assume Conjecture \ref{con:xz-conjecture_SO(N)} is true. Then Mathieu's conjecture is true for $SO(N)$.
\end{proposition}
\begin{proof}[Outline of the proof]
	The proof goes analogously to the proof we gave to prove Theorem \ref{thm:Mathieu_proven_assuming_XZ-conjecture}, with some simplifications. Let $f,h$ be finite-type functions of $SO(N)$, and assume $\int_{SO(N)} f^P(g) dg=0$ for all $P\in \N$. Then by applying Lemma \ref{lemma:Main_Lemma_SO(N)} and Conjecture \ref{con:xz-conjecture_SO(N)} we see that $0$ does not lie in the convex hull of $\mathrm{Sp}(\tilde{f})$. Now assume that $\int_{SO(N)} f^Ph \neq 0$ for infinitely many $P\in \N$. Because of linearity, we can assume that $h$ is a monomial with respect to each component. In other words
	\begin{align*}\begin{split}
	f(g)=f_{SO(N)}(g)&= \sum_{j=1}^M\sum_{i=1}^Q c_{ij} e^{ik_{ij}^1\phi_1}\cos^{m_{ij}^1}(\phi_2)\sin^{n_{ij}^1}(\phi_2)\cdots \cos^{m_{ij}^{N-1}}(\phi_{N-1})\sin^{n_{ij}^{N-1}}(\phi_{N-1})\\
	&\cdot (f_{SO(N-1)})_{ij}(g_{SO(N-1)})\end{split}\\
	\begin{split}h(g)=h_{SO(N)}(g)= C& e^{iK_1\phi_1}\cos^{R_1}(\phi_2)\sin^{S_1}(\phi_2)\cdots \cos^{R_{N-1}}(\phi_{N-1})\sin^{S_{N-1}}(\phi_{N-1}) (h_{SO(N-1)})(g_{SO(N-1)}),\end{split}
	\end{align*}	
	where $h_{SO(N-1)}$ is a monomial finite-type function on $SO(N-1)$. If $\int f^Ph\neq 0$ then there exists at least one set of parameters $\{\beta_{ij}\}_{i,j}\in\N_0$ where $1\leq i\leq M$ and $1\leq j\leq Q$ with $\sum_{i,j}\beta_{ij}=P$ such that
	\begin{align*}
	0\neq& \int_{[-1,1]^{N-2}}\int_{(S^1)^{N-1}} e^{i(\sum_{i,j}\beta_{ij}k_{ij}^1+K_1)\phi_1}\cos^{\sum_{i,j}\beta_{ij}m_{ij}^1+R_1}(\phi_2)\sin^{\sum_{i,j}\beta_{ij}n_{ij}+S_1}(\phi_2)\cdots\\
	&\cdot \cos^{\sum_{i,j}\beta_{ij}m_{ij}^{N-1}+R_{N-1}}(\phi_{N-1})\sin^{\sum_{i,j}\beta_{ij}n_{ij}+S_{N-1}}(\phi_{N-1}) \left(\prod_{j=1}^{N-1}\sin^{j-1}(\phi_j)\right)d\phi_1\ldots d\phi_{N-1}\\
	&\cdot\left(\prod_{i,j}\int_{SO(N-1)} [f_{SO(N-1)}(g_{SO(N-1)})]^{\beta_{ij}} h_{SO(N-1)}(g_{SO(N-1)})dg_{SO(N-1)}\right).
	\end{align*}
	Note that $dg$ is a Haar measure, hence $\int_G f(g) dg = \int_G f(xg) dg$ for any $x\in G$. Sending $g\mapsto xg$ where $x=e^{\psi\lambda_2}$ for any $\psi\in[0,2\pi)$ gives, taking Lemma \ref{lemma:Euler_Angles_SO(N)} into account, that this is equivalent to replacing $\phi_1$ with $\phi_1+\psi$. Applying this to the above integral, stating that it is invariant under this mapping, we have 
	$$e^{i(\sum_{i,j}\beta_{ij}k_{ij}^1+K_1)\psi}=1$$ for all $\psi\in[0,2\pi)$. This can only be fulfilled if $$\sum_{i,j}\beta_{ij}k_{ij}^1+K_1=0.$$ Since $dg_{SO(N-1)}$ is also a Haar measure, applying inductively the same strategy to the integral $$\sum_{i,j}\int_{SO(N-1)}[f_{SO(N-1)}]^{\beta_{ij}}h_{SO(N-1)} dg_{SO(N-1)}$$ gives $\sum_{i,j}\beta_{ij}k_{ij}^q+K_q=0$ for all $1\leq q\leq N-1$. But this means that $$-(K_1,\ldots,K_{N-1})=\sum_{ij}\beta_{ij}(k_{ij}^1,\ldots,k_{ij}^{N-1}).$$ Dividing both sides by $P$ and taking the limit $P\rightarrow\infty$, we see that the left-hand side goes to 0. Since the convex hull is closed, we must have that $0$ lies in the convex hull of $\mathrm{Sp}(\tilde{f})$ which is a contradiction. Hence $\int_{SO(N)} f^Ph=0$ for all $P$ large enough, hence proving Mathieu's conjecture.
	
\end{proof}

	\newpage
	\appendix 

\section{Appendix}
In this appendix, we will prove Lemma \ref{lemma:Euler_Angles} and Lemma \ref{lemma:Haar_Measure}. We will note that the proofs are inspired by and/or based on \cite{bertini2006euler,Parametrization_SU(4),Parametrization_SU(N),spengler2012composite}.

\begin{proof}[Proof of Lemma \ref{lemma:Euler_Angles}]
	Since the mapping in Lemma \ref{lemma:Euler_Angles} is defined inductively, the proof will be by induction on $N$. We start with $N=2$. On $SU(2)$ it can easily be seen that the function $$F_2(\phi_1,\psi_1,\omega_1)=\begin{pmatrix}
	e^{i(\phi_1+\omega_1)}\cos(\psi_1)&e^{i(\phi_1-\omega_1)}\sin(\psi_1)\\
	e^{-i(\phi_1-\omega_1)}\sin(\psi_1)&e^{-i(\phi_1+\omega_1)}\cos(\psi_1)
	\end{pmatrix}$$ is surjective, and a bijection between $[0,\pi]\times\left[0,\frac{\pi}{2}\right]\times[0,2\pi]$ and $SU(2)$ up to measure zero sets. So let it be proven for $N-1$. Let us consider $N$. We note that $\left(\prod_{2\leq k\leq {N-1}}A(k)(\phi_k,\psi_k)\right)\in \begin{pmatrix}
	SU(N-1)&0\\
	0&1
	\end{pmatrix}$. For simplicity, let us denote the following two matrices 
	\begin{align}
	\begin{pmatrix}
	A&0\\
	0&1
	\end{pmatrix}:=\left(\prod_{2\leq k\leq N-1}A(k)(\phi_k,\psi_k)\right)e^{\lambda_3\phi_{N-1}},\quad B:=F_{N-1}(\phi_{N},\ldots,\phi_{\frac{N(N-1)}{2}},\psi_{N},\ldots,\psi_{\frac{N(N-1)}{2}},\omega_1,\ldots,\omega_{N-2}).
	\end{align} Then we can write $F_{N}(\phi_1,\ldots,\omega_{N-1})$ as  
	\begin{align*}
	F_{N}(\phi_{1},\ldots,\omega_{N-1})=\begin{pmatrix}
	A&0\\
	0&1
	\end{pmatrix}\begin{pmatrix}
	\cos(\psi_{N-1})&0&\ldots&0&\sin(\psi_{N-1})\\
	0&1&&&0\\
	\vdots&&\ddots&&\vdots\\
	0&&&1&0\\
	-\sin(\psi_{N-1})&0&\ldots&0&\cos(\psi_{N-1})
	\end{pmatrix}\begin{pmatrix}
	B&0\\
	0&1
	\end{pmatrix}\begin{pmatrix}
	e^{i\omega_{N-1}}\cdot\mathbf{1}_{N-1}&\\
	&e^{-i(N-1)\omega_{N-1}}
	\end{pmatrix}.
	\end{align*}
	We will start with proving that $F_N$ is injective on the interior of the hypercube. Let $X:=F_{N}(\phi_1,\ldots,\omega_{N-1})$ and $Y:=F_{N}(\phi_1',\ldots,\omega_{N-1}')$ satisfy $X=Y$. Since $X,Y\in SU(N)$ we have $XY^{-1}=X^{-1}Y=\mathbf{1}_N$. If we look at the matrix component $(XY^{-1})_{N,N}$, we get the following equation:
	\begin{align}\label{eq:Injectivity_Euler_parametrization}
	\begin{split}e^{i(\omega_{N-1}-\omega_{N-1}')}\cdot\left[\begin{pmatrix}
	-\sin(\psi_{N-1})&0&\ldots&0
	\end{pmatrix}B(B')^{-1}\begin{pmatrix}
	-\sin(\psi_{N-1}')\\
	0\\
	\vdots\\
	0
	\end{pmatrix}\right]+&\\ \cos(\psi_{N-1})\cos(\psi_{N-1}')e^{-i(N-1)(\omega_{N-1}-\omega_{N-1}')}&=1.\end{split}
	\end{align}
	Taking the absolute value on both sides, and using the triangle inequality, we get
	\begin{align*}
	|\sin(\psi_{N-1})\sin(\psi_{N-1}')|\cdot\left|\left(B(B')^{-1}\right)_{11}\right|+|\cos(\psi_{N-1})\cos(\psi_{N-1}')|\geq 1.
	\end{align*}
	Since we know that $B(B')^{-1}\in SU({N-1})$, we have that each matrix element has $|\left(B(B')^{-1}\right)_{ij}|\leq 1$. In addition, $\psi_{N-1},\psi_{N-1}'\in [0,\frac{\pi}{2}]$ thus $\cos(\psi_{N-1}),\cos(\psi_{N-1}'),\sin(\psi_{N-1}),\sin(\psi_{N-1}')\geq 0$. Therefore we find 
	\begin{align*}
	1&\leq |\sin(\psi_{N-1})\sin(\psi_{N-1}')|\cdot\left|\left(B(B')^{-1}\right)_{11}\right|+|\cos(\psi_{N-1})\cos(\psi_{N-1}')|\\
	&\leq \sin(\psi_{N-1})\sin(\psi_{N-1}')+\cos(\psi_{N-1})\cos(\psi_{N-1}')\\
	&=\cos(\psi_{N-1}-\psi_{N-1}')\leq 1.
	\end{align*}
	Hence we see that all inequalities were equalities. So $\cos(\psi_{N-1}-\psi_{N-1}')=1$, which can only mean $\psi_{N-1}=\psi_{N-1}'$ in the given range. In addition, it means $|\left(B(B')^{-1}\right)_{11}|=1$. So set $\left(B(B')^{-1}\right)_{11}=e^{i\xi}$ for some $\xi\in\R$. Putting it all in Equation (\ref{eq:Injectivity_Euler_parametrization}) and writing the real and imaginary part separately gives
	\begin{align}\label{eq:Injectivity_Euler_parametrization2}
	\begin{split}
	\cos(\omega_{N-1}-\omega_{N-1}'+\xi)\sin^2(\psi_{N-1})+\cos((N-1)(\omega_{N-1}-\omega_{N-1}'))\cos^2(\psi_{N-1})&=1,\\
	\sin(\omega_{N-1}-\omega_{N-1}'+\xi)\sin^2(\psi_{N-1})-\sin((N-1)(\omega_{N-1}-\omega_{N-1}'))\cos^2(\psi_{N-1})&=0.
	\end{split}
	\end{align}
	Taking the absolute value of the first equation, and again using the triangle inquality gives
	\begin{align*}
	1&\leq |\cos(\omega_{N-1}-\omega_{N-1}'+\xi)|\sin^2(\psi_{N-1})+|\cos((N-1)(\omega_{N-1}-\omega_{N-1}'))|\cos^2(\psi_{N-1})\\
	&\leq \sin^2(\psi_{N-1})+\cos^2(\psi_{N-1})=1.
	\end{align*}
	Once again, we see that $|\cos((N-1)(\omega_{N-1}-\omega_{N-1}'))|=|\cos(\omega_{N-1}-\omega_{N-1}'+\xi)|=1$. But since $\sin^2,\cos^2\geq 0$, we must have that $\cos((N-1)(\omega_{N-1}-\omega_{N-1}'))=\cos(\omega_{N-1}-\omega_{N-1}'+\xi)=1$ in order to satisfy Equation (\ref{eq:Injectivity_Euler_parametrization2}). We have two solutions. Either $\omega_{N-1}=\omega_{N-1}'$ or $\omega_{N-1}=\frac{2\pi}{N-1}$ and $\omega_{N-1}'=0$ or vice versa. However, since we are considering the interior of the hypercube, the latter cannot be the case. Therefore, we get $\omega_{N-1}=\omega_{N-1}'$, and so $\xi=0$. Going back to $XY^{-1}=\mathbf{1}_N$, using what we found, gives	
	\begin{align*}
	\mathbf{1}_N=\begin{pmatrix}
	A&0\\
	0&1
	\end{pmatrix}&\begin{pmatrix}
	\cos(\psi_{N-1})&0&\sin(\psi_{N-1})\\
	0&\mathbf{1}_{N-2}&0\\
	-\sin(\psi_{N-1})&0&\cos(\psi_{N-1})
	\end{pmatrix}\begin{pmatrix}
	B(B')^{-1}&0\\
	0&1
	\end{pmatrix}\\
	&\cdot\begin{pmatrix}
	\cos(\psi_{N-1})&0&-\sin(\psi_{N-1})\\
	0&\mathbf{1}_{N-2}&0\\
	\sin(\psi_{N-1})&0&\cos(\psi_{N-1})
	\end{pmatrix}\begin{pmatrix}
	(A')^{-1}&\\
	&1
	\end{pmatrix}.
	\end{align*}
	Note that this can be rewritten as $P\begin{pmatrix}
	B(B')^{-1}&0\\
	0&1
	\end{pmatrix}P^{-1}=\begin{pmatrix}
	A^{-1}A'&0\\
	0&1
	\end{pmatrix}$ where $P\in SO(N)$ is a rotation in the $(e_1,e_{N})$-plane. This means that either $P=1$ or $(B(B')^{-1})_{1k}=(B(B')^{-1})_{k1}=\delta_{1k}$ for all $k=1,\ldots,N-1$. Either case immediately gives $AB(B')^{-1}(A')^{-1}=1$. But that would mean $(B')^{-1}(A')^{-1}=(A'B')^{-1}$ is the inverse of $AB$, so $AB=A'B'$.
	
	By the same arguments, considering $X^{-1}Y$, will give $BA=B'A'$. Note that $P=1$ can only occur when $\psi_{N-1}=0$, which would lie at the boundary of the hypercube, which we are not considering. So we can assume $(B(B')^{-1})_{1k}=(B(B')^{-1})_{k1}=\delta_{1k}$ for all $k$. Since $AB(B')^{-1}(A')^{-1}=\mathbf{1}_N$, it gives $(A(A')^{-1})_{1k}=(A(A')^{-1})_{k1}=\delta_{1k}$ as well. One can do the calculations and see that this is only possible if $\phi_{N-1}=\phi_{N-1}'$ and $\phi_j=\phi_j'$ and $\psi_j=\psi_j'$ for all $j=1,\ldots N-2$. In other words $A=A'$ and thus $B=B'$. This means $$F_{N-1}(\phi_N\ldots,\omega_{N-2})=B=B'=F_{N-1}(\phi_N',\ldots,\omega_{N-2}').$$ Using the induction hypothesis we see that $\phi_j=\phi_j'$, $\psi_j=\psi_j'$ for all $j=N,\ldots,\frac{N(N-1)}{2}$ and $\omega_k=\omega_k'$ for all $k=1,\ldots,N-2$. Therefore, we see that $F_{N}$ is injective on the interior of the hypercube, and not injective on the boundary.	
	
	Next, we consider surjectivity. Let $U\in SU(N)$. If we can show that there exist $\phi_1,\ldots,\omega_{N-1}$ such that $F_N(\phi_1,\ldots,\omega_{N-1})^{-1}U=\mathbf{1}_{N}$, then $U=F_N(\phi_1,\ldots,\omega_{N-1})$. The last three matrices of the product $F_N(\phi_1,\ldots,\omega_{N-1})^{-1}U$ are of the form
	\begin{align*}
	\begin{pmatrix}
	\cos(\psi_{1})&-\sin(\psi_1)&&&\\
	\sin(\psi_1)&\cos(\psi_1)&&&\\
	&&1&\\
	&&&\ddots&&\\
	&&&&1
	\end{pmatrix}\begin{pmatrix}
	e^{-i\phi_1}&&&&\\
	&e^{i\phi_1}&&&\\
	&&1&&\\
	&&&\ddots&\\
	&&&&1
	\end{pmatrix}\begin{pmatrix}
	u_{11}&u_{12}&\ldots &u_{1N}\\
	u_{21}&u_{22}&\ldots &u_{2N}\\
	\vdots&\vdots&\ddots &\vdots\\
	u_{N1}&u_{N2}&\ldots& u_{NN}
	\end{pmatrix},
	\end{align*}
	where $U$ is written as $U=(u_{ij})_{i,j=1,\ldots, N}$. Multiplying these three matrices gives a new element in $SU(N)$, call this $U'$. Doing the multiplication gives
	\begin{align*}
	u_{1N}'&=u_{1N}e^{-i\phi_1}\cos(\psi_1)-u_{2N}e^{i\phi_1}\sin(\psi_1),\\
	u_{2N}'&=u_{1N}e^{-i\phi_1}\sin(\psi_1)+u_{2N}e^{i\phi_1}\cos(\psi_1).
	\end{align*}
	We choose $\psi_1,\phi_1$ in such a way that $u_{2N}'=0$. This is always possible. For if $u_{1N}=0$ then we choose $\psi_1=\frac{\pi}{2}$ and $\phi_1=0$. If $u_{2,N}=0$ we choose $\psi_1=0$ and $\phi_1=0$. If both are equal to zero, we are free to choose what we want, so we choose $\psi_1=\phi_1=0$. And if neither is 0, we choose $\tan(\psi_1)=\left|\frac{u_{2N}}{u_{1N}}\right|$ and $2\phi_1=\arg\left(\frac{-u_{1N}}{u_{2N}}\right)$. Since $\phi_1\in[0,\pi)$, these determine $\phi_1,\psi_1$ uniquely.
	
	Filling in our choice of $\psi_1,\phi_1$ shows that the next two matrices in the multiplication $F_N(\phi_1,\ldots,\omega_{N-1})^{-1}U$ will be of the form 
	\begin{align*}
	\begin{pmatrix}
	\cos(\psi_{2})&0&-\sin(\psi_2)&&&\\
	0&1&0&&&\\
	\sin(\psi_2)&0&\cos(\psi_2)&&&\\
	&&&1&\\
	&&&&\ddots&&\\
	&&&&&1
	\end{pmatrix}\begin{pmatrix}
	e^{-i\phi_2}&&&&\\
	&e^{i\phi_2}&&&\\
	&&1&&\\
	&&&\ddots&\\
	&&&&1
	\end{pmatrix}\begin{pmatrix}
	u_{11}'&u_{12}'&\ldots &u'_{1N}\\
	u_{21}'&u_{22}'&\ldots &0\\
	u_{31}'&u_{32}'&\ldots &u'_{3N}\\
	\vdots&\vdots&\ddots &\vdots\\
	u_{N1}'&u_{N2}'&\ldots& u_{NN}'
	\end{pmatrix}.
	\end{align*}
	If we denote this matrix product by $U''=(u''_{ij})_{i,j=1,\ldots,N}$ we get
	\begin{align*}
	u_{1N}''&=u_{1N}'e^{-i\phi_2}\cos(\psi_2)-u_{3N}'\sin(\psi_2),\\
	u_{3N}''&=u_{1N}'e^{-i\phi_2}\sin(\psi_2)+u_{3N}'\cos(\psi_2).
	\end{align*}
	To set $u_{3N}''=0$, we consider a few cases. If $u_{1N}'=0$ we set $\psi_{2}=\frac{\pi}{2}$ and $\phi_2=0$. If $u_{3N}'=0$ then $\psi_2=0$ and $\phi_2=0$. And if neither are 0, we set $\tan(\psi_2)=\left|\frac{u'_{3N}}{u'_{1N}}\right|$ and $\phi_3=\arg\left(\frac{-u_{1N}}{u_{3N}}\right)$. Since $\psi_2\in\left[0,\frac{\pi}{2}\right]$ and $\phi_2\in[0,2\pi]$, they are determined uniquely in the last case.
	
	This procedure can be repeated $N-3$ more times, defining $\psi_1,\ldots,\psi_{N-2},\phi_1,\ldots,\phi_{N-2}$, and denoting the resulting matrix by $U^{(N-2)}=(u_{ij}^{(N-2)})_{i,j=1,\ldots,N}$. Then the next set of multiplications is given by
	\begin{align*}
	\begin{pmatrix}
	\cos(\psi_{N-1})&0&-\sin(\psi_{N-1})\\
	0&\mathbf{1}_{N-2}&0\\
	\sin(\psi_{N-1})&0&\cos(\psi_{N-1})
	\end{pmatrix}\begin{pmatrix}
	e^{-i\phi_{N-1}}&&\\
	&e^{i\phi_{N-1}}&\\
	&&\mathbf{1}_{N-2}
	\end{pmatrix}\begin{pmatrix}
	u_{11}^{(N-2)}&u_{12}^{(N-2)}&\ldots&u^{(N-2)}_{1(N-1)}&u^{(N-2)}_{1N}\\
	u_{21}^{(N-2)}&u_{22}^{(N-2)}&\ldots&u^{(N-2)}_{1(N-1)}&0\\
	\vdots&\vdots&\ddots&\vdots&\vdots\\
	u_{(N-1)1}^{(N-2)}&u_{(N-1)2}^{(N-2)}&\ldots&u^{(N-2)}_{1(N-1)}&0\\
	u_{N1}^{(N-2)}&u_{N2}^{(N-2)}&\ldots&u^{(N-2)}_{N(N-1)}&u^{(N-2)}_{NN}\\
	\end{pmatrix}.
	\end{align*}
	Denoting the resulting matrix by $U^{(N-1)}=(u_{ij}^{N-1})_{i,j=1,\ldots,N}$ gives the following set of equations:
	\begin{align*}
	u^{(N-1)}_{1N}&=u^{(N-2)}_{1N}e^{-i\phi_{N-1}}\cos(\psi_{N-1})-u^{(N-2)}_{NN}\sin(\psi_{N-1}),\\
	u^{(N-1)}_{NN}&=u^{(N-2)}_{1N}e^{-i\phi_{N-1}}\sin(\psi_{N-1})+u^{(N-2)}_{NN}\cos(\psi_{N-1}).
	\end{align*}
	We now choose $u^{(N-1)}_{1N}=0$. This can be achieved in a similar way as before. Therefore the resulting matrix is then
	\begin{align*}
	U^{(N-1)}=\begin{pmatrix}
	u_{11}^{(N-1)}&\ldots &u^{(N-1)}_{1(N-1)}&0\\
	u_{21}^{(N-1)}&\ldots &u^{(N-1)}_{2(N-1)}&0\\
	\vdots&\vdots &\vdots&\vdots\\
	u_{(N-1)1}^{(N-1)}&\ldots &u^{(N-1)}_{(N-1)(N-1)}&0\\
	u_{N1}^{(N-1)}&\ldots &u^{(N-1)}_{N(N-1)}&u^{(N-1)}_{NN}
	\end{pmatrix}.
	\end{align*}
	Since $U^{(N-1)}\in SU(N)$, we have $U^{(N-1)}(U^{(N-1)})^\dagger =\mathbf{1}_N$. If we write $X=(u^{(N-1)}_{ij})_{i,j=1,\ldots,N-1}$ we get 
	\begin{align*}
	\begin{pmatrix}
	X&\begin{pmatrix}
	0\\
	\vdots\\
	0
	\end{pmatrix}\\
	\begin{pmatrix}
	u^{(N-1)}_{N1}&\ldots&u^{(N-1)}_{N(N-1)}
	\end{pmatrix}& u^{(N-1)}_{NN}
	\end{pmatrix}
	\begin{pmatrix}
	X^\dagger&\begin{pmatrix}
	\overline{u^{(N-1)}_{N1}}\\
	\vdots\\
	\overline{u^{(N-1)}_{N(N-1)}}
	\end{pmatrix}\\
	\begin{pmatrix}
	0&\ldots&0
	\end{pmatrix}& \overline{u^{(N-1)}_{NN}}
	\end{pmatrix}=\mathbf{1}_N.
	\end{align*}
	This implies $XX^\dagger=\mathbf{1}_{N-1}$, which means by finite dimensionality arguments that $X\in U(N-1)$. We also see that 
	$$X\begin{pmatrix}
	\overline{u^{(N-1)}_{N1}}\\
	\vdots\\
	\overline{u^{(N-1)}_{N(N-1)}}
	\end{pmatrix}=0.$$ Since $X$ is invertible, we have $u^{(N-1)}_{Nk}=0$ for $k=1,2,\ldots,N-1$. Therefore $$U^{(N-1)}=\begin{pmatrix}
	X&0\\
	0&u^{(N-1)}_{NN}
	\end{pmatrix}.$$ 
	We recall we were originally looking at $[F_{N}(\phi_1,\ldots,\omega_{N-1})]^{-1}U=\mathbf{1}_N$. Applying our procedure, the multiplication is reduced to finding the remaining parameters such that
	\begin{align*}
	\mathbf{1}_N=\begin{pmatrix}
	e^{-i\omega_{N-1}}&&&\\
	&\ddots&&\\
	&&e^{-i\omega_{N-1}}&\\
	&&&e^{i(N-1)\omega_{N-1}}
	\end{pmatrix}
	\begin{pmatrix}
	[F_{N-1}(\phi_N,\ldots,\omega_{N-2})]^{-1}&0\\
	0&1
	\end{pmatrix}
	\begin{pmatrix}
	X&\\
	&u^{(N-1)}_{NN}
	\end{pmatrix}.
	\end{align*}
	We note that the two left-most matrices commute, hence this is equivalent to
	\begin{align*}
	\mathbf{1}_N=\begin{pmatrix}
	[F_{N-1}(\phi_N,\ldots,\omega_{N-2})]^{-1}&0\\
	0&1
	\end{pmatrix}
	\begin{pmatrix}
	e^{-i\omega_{N-1}}&&&\\
	&\ddots&&\\
	&&e^{-i\omega_{N-1}}&\\
	&&&e^{i(N-1)\omega_{N-1}}
	\end{pmatrix}\begin{pmatrix}
	X&\\
	&u^{(N-1)}_{NN}
	\end{pmatrix}.
	\end{align*}
	Since $X\in U(N-1)$ we have $\det(X)=e^{i\xi}$ for some $\xi\in[0,2\pi)$. In addition, we see that $$1=\det(U^{(N-1)})=\det(X)u_{NN}^{(N-1)}=e^{i\xi}u_{NN}^{(N-1)},$$ so $u_{NN}^{(N-1)}=e^{-i\xi}$. Choosing $\omega_N=\frac{\xi}{N-1}$ gives us the equation
	\begin{align*}
	\mathbf{1}_N=\begin{pmatrix}
	F_N(\phi_N,\ldots,\omega_{N-2})&0\\
	0&1
	\end{pmatrix}
	\begin{pmatrix}
	e^{\frac{-i\xi}{N-1}}X&\\
	&1
	\end{pmatrix}.
	\end{align*}
	By the induction hypothesis, we know $F_{N-1}$ is surjective onto $SU(N-1)$. Since $\det(e^{\frac{-i\xi}{N-1}}X)=e^{-i\xi}e^{i\xi}=1$, we get that $e^{\frac{-i\xi}{N-1}}X\in SU(N-1)$ hence it can be reached by $F_{N-1}$. Therefore, $F_{N}$ is surjective.
	
	Finally we want to show that $F_N$ is a diffeomorphism. It is clear that $F_N$ is $C^\infty$, because the exponential and multiplication are $C^\infty$. On the interiour of the hypercube, we see that the inverse is given inductively by $$\psi_k=\arctan\left|\frac{u^{(k-1)}_{k+1,N}}{u^{(k-1)}_{1,N}}\right|,\qquad \phi_k=\frac{1}{2^{\delta_{k,1}}}\arg\left(\frac{-u^{(k-1)}_{1,N}}{u^{(k-1)}_{k-1,N}}\right),\qquad\omega_N=\frac{\arg(\det((u^{(N-1)}_{ij})_{i,j=1,\ldots,N-1}))}{N-1}.$$ We see that each of these equations are continuous and differentiable on the image of the interior of the domain of $F_N$. In other words, $F_N$ is on the interior of the hypercube a diffeomorphism.
\end{proof}


\begin{proof}[Proof of Lemma \ref{lemma:Haar_Measure}]
	The proof is based on \cite{bertini2006euler}. Before proving the lemma in detail, we outline the general strategy. We will construct the proof in four steps. First, we will show that there is a closed subgroup $K$ such that $G/K$ is a symmetric space. We wish to show $$\int_G f(g)dg = \int_{G/K}\int_K f(xk)dk\,dg_K$$ for any measurable function $f$ on $G$, and where $x\in G$ is a representative of $xK\in G/K$ and $k\in K$ such that $g=xk$. Here $dk$ is the Haar measure on $K$ and $dg_K$ is the unique $G$-invariant measure on the symmetric space $G/K$ \cite{Helgason_DifGeom_Symm_Spaces}. If the previous equation is true, then it shows that $$dg=dk\,dg_K.$$
	
	Second, we construct left-invariant one-forms on $G/K$, which can be wedged to find $dg_K$. Third, we will show how the top form $dg_K$ looks like explicitly by considering the parameterization of $SU(N)$ as in Lemma \ref{lemma:Euler_Angles}. We end the proof by normalizing the measure to get the Haar measure, which we shall call $dg_{SU(N)}.$
	
	Our first step will be to find the subgroup $K$. Let us consider the group $$K:=\left\{\left.\begin{pmatrix}
	A&0\\
	0&1
	\end{pmatrix}\begin{pmatrix}
	e^{i\omega_{N}}\mathbf{1}_{N-1}&0\\
	0&e^{-i(N-1)\omega_{N-1}}
	\end{pmatrix}\,\right|\,A\in SU(N-1), \omega_{N-1}\in \left[0,\frac{2\pi}{N-1}\right]\right\}\simeq U(N-1).$$ Note that this subgroup is closed, and is the same subgroup as for the $KAK$ decomposition, as discussed in Remark \ref{remark:KAK_decomposition}. This automatically shows that $(G,K)$ is a Riemannian symmetric pair and thus $G/K$ a symmetric space. To show the identity $$\int_G f(g)dg = \int_{G/K}\int_K f(xk)dk\,dg_K,$$ where $x\in G$ is a representative of $xK\in G/K$ and $k\in K$, it is enough to show that $|\det(\Ad_G(k))|=|\det(\Ad_K(k))|$ \cite{Helgason_DifGeom_Symm_Spaces}. Since we are considering matrix Lie groups, $\Ad_G(k)(X)=kXk^{-1}$ for any $X\in\mathfrak{g}$. The Lie algebra of $K$, which we denoted as $\mathfrak{k}$, is generated by $\lambda_1,\lambda_2,\ldots,\lambda_{(N-1)^2-1},\lambda_{N^2-1}.$ Let us denote $\mathfrak{p}:=\Span_\R(\lambda_{(N-1)^2},\ldots,\lambda_{N^2-2})$. We see then that for any $k\in K$ and $1\leq l\leq 2(N-1)$ that
	\begin{align}\label{eq:Proof_Haar_measure_Adjoint_lies_in_SO(2(N-1))}
	k\lambda_{(N-1)^2-1+l}k^{-1}&=\begin{pmatrix}
	A&0\\
	0&1
	\end{pmatrix}\begin{pmatrix}
	e^{i\omega_{N-1}}\mathbf{1}_{N-1}&0\\
	0&e^{-i(N-1)\omega_{N-1}}
	\end{pmatrix}\begin{pmatrix}
	\mathbf{0}&\vec{v}\\
	-(\vec{v})^\dagger&0
	\end{pmatrix}\begin{pmatrix}
	e^{-i\omega_{N-1}}\mathbf{1}_{N-1}&0\\
	0&e^{i(N-1)\omega_{N-1}}
	\end{pmatrix}\begin{pmatrix}
	A^{-1}&0\\
	0&1
	\end{pmatrix}\nonumber\\
	&=\begin{pmatrix}
	A&0\\
	0&1
	\end{pmatrix}\begin{pmatrix}
	\mathbf{0}&e^{iN\omega_{N-1}}\vec{v}\\
	-e^{-iN\omega_{N-1}}(\vec{v})^\dagger&0
	\end{pmatrix}\begin{pmatrix}
	A^{-1}&0\\
	0&1
	\end{pmatrix}\nonumber\\
	&=\begin{pmatrix}
	\mathbf{0}&Ae^{iN\omega_{N-1}}\vec{v}\\
	-(\vec{v})^\dagger e^{-iN\omega_{N-1}}A^{-1}&0
	\end{pmatrix} = \begin{pmatrix}
	\mathbf{0}&Ae^{iN\omega_{N-1}}\vec{v}\\
	-(Ae^{iN\omega_{N-1}}\vec{v})^\dagger&0
	\end{pmatrix}.
	\end{align}
	where $\lambda_{(N-1)^2-1+l}=\begin{pmatrix}
	\mathbf{0}&\vec{v}\\
	-(\vec{v})^\dagger&0
	\end{pmatrix},$ $\vec{v}\in \C^{N-1}$ is a column vector, and $(\vec{v})^\dagger$ is the adjoint of $\vec{v}$ i.e. the transpose of the complex conjugate of $\vec{v}$. Therefore $\Ad_G(k)\lambda_{(N-1)^2-1+l}\in \mathfrak{p}$ for all $l$ and thus
	\begin{align}\label{eq:Proof_Haar_Ad(h)_acts_diagonally}
	\Ad_G(k)=\begin{pmatrix}
	\Ad_K(k)&0\\
	0&\Ad_G(k)|_{\mathfrak{p}}
	\end{pmatrix}.
	\end{align}
	Note that we can identify $\mathfrak{p}\simeq \R^{2(N-1)}$. Also recall that $U(N)\subset SO(2N)$ by identifying $\C^n\simeq \R^{2n}$. By looking at Equation (\ref{eq:Proof_Haar_measure_Adjoint_lies_in_SO(2(N-1))}) it can be seen that $\Ad_G(k)|_{\mathfrak{p}}\subset SO(2(N-1))$. Hence we conclude that
	$$|\det(\Ad_G(k))|=|\det(\Ad_K(k))|\,|\det(\Ad_G(k)|_{\mathfrak{p}})|=|\det(\Ad_K(k))|.$$
	Therefore there exists a unique $G$-invariant measure on $G/K$ and 
	\begin{align*}
	dg=dg_K\,dk=dg_{K}\;dg_{SU(N-1)}\,d\omega_{N-1}.
	\end{align*}
	The rest of the proof will be dedicated to finding $dg_K$. To find $dg_{K}$ explicitly, we will consider the Maurer-Cartan one-form $\omega$, which, at $g\in G$, is a map $\omega_g:T_gG\rightarrow \mathfrak{g}$ defined by
	$$\omega_g(X_g):=T_g(L_{g^{-1}})X_g.$$
	Note that by construction $\omega_g$ is left-invariant, that is to say $(L_x)^*\omega_g=\omega_g$ for all $g\in G$. In the case of matrix groups, especially when $G=SU(N)$, $\omega_g$ can be calculated explicitly to be 
	\begin{align*}
	\omega_g=g^{-1}dg=\sum_{j=1}^{N^2-1}g^{-1}\frac{\partial g}{\partial x_j}dx_j,
	\end{align*} where $x_1,\ldots,x_{N^2-1}$ is a set of local coordinates of $G$.
	We recall that $\Tr(\lambda_j\lambda_k)=0$ whenever $j\neq k$. Using this, we construct one-forms out of $\omega_g$ by defining for $m=1,\ldots,2(N-1)$ the form
	\begin{align*}
	(e^m)_g:T_gG\rightarrow \R,\qquad (e^m)_g(X_g):=\frac{1}{2}\Tr(\omega_g(X_g)\lambda_{(N-1)^2-1+m}).
	\end{align*}
	Note that $e^m$ is left-invariant, because $\omega_g$ is left-invariant. Let $g=xk$ where $k\in K$ and $x\in G$. Then we see that 
	$$\omega_g=(xk)^{-1}d(xk)=k^{-1}x^{-1}(dx)\,k+k^{-1}dk.$$
	Filling this into $e^m$, and noting that $\Tr(w\lambda_{(N-1)^2-1+m})=0$ for all $w\in \mathfrak{k}$, gives
	\begin{align*}
	(e^m)_{xk}=\frac{1}{2}\Tr\left(k^{-1}x^{-1}dx\,k\lambda_{(N-1)^2-1+m}\right)= \frac{1}{2}\Tr\left(x^{-1}dx\,\Ad_G(k)(\lambda_{(N-1)^2-1+m})\right).
	\end{align*}
	Now let us define the $2(N-1)$-form given by
	\begin{align*}
	\mu:=e^1\wedge\ldots\wedge e^{2(N-1)}.
	\end{align*}
	This form is left-invariant because $e^m$ is left-invariant. In addition, if $k\in K$, then by Equation (\ref{eq:Proof_Haar_Ad(h)_acts_diagonally}) we see  $$\mu_{xk}=\det(\Ad_G(k)|_\mathfrak{p})\,\mu_{x1}.$$ Note that $\det(\Ad_G(k)|_{\mathfrak{p}})\in S^1$, but $\mu_g$ has values in $\R$ for all $g\in G$. Therefore in order for the above equation to hold, one must have $\det(\Ad_G(k)|_{\mathfrak{p}})=\pm 1$ for all $k\in K$. But $K$ is connected and $\det$ and $\Ad_K$ are smooth mappings, hence the image is a connected set. Since $\Ad_G(1)=\text{Id}$, we must have $\det(\Ad_G(k)|_\mathfrak{p})=1$. Therefore we conclude that $$\mu_{xk}=\mu_{x1}.$$ Together with the fact that $\mathfrak{k}$ lies in the kernel of $(e^m)_e$, the form $\mu$ can be identified as a top form on the symmetric space $G/K$. In addition, $\mu$ is $G$-invariant, hence we can conclude that 
	\begin{align*}
	(dg_K)_{xK}=\mu_{x1} = (e^1)_{x1}\wedge\ldots\wedge (e^{2N-1})_{x1}.
	\end{align*}  where
	\begin{align} \label{eq:Proof_Haar_zo_ziet_e^k_er_in_praktijk_uit}
	(e^m)_{x1}=\frac{1}{2}\Tr(x^{-1}dx\lambda_{(N-1)^2-1+m}).
	\end{align}
	If $x_1,\ldots,x_{2(N-1)}$ is a parametrization of $G/K$, we have in general that $$e^1\wedge\ldots\wedge e^{2(N-1)}=\det((e_{ij}))\,dx_1\wedge\ldots\wedge dx_{2(N-1)}$$
	where $e^m=\sum_{j=1}^{2(N-1)}e_{mj}dx_{j}$. By Lemma \ref{lemma:Euler_Angles}, we have that the set $\phi_{N-1},\psi_{N-1},\ldots,\phi_1,\psi_1$ parametrize $G/K$. Hence we get that
	\begin{align}\label{eq:Proof_Haar_measure_equation_that_only_matters}
	dg_K &= e^1\wedge \ldots\wedge e^{2N}\nonumber\\
	&= \det((e_{ij}))\,d\phi_{N-1}\wedge d\psi_{N-1}\wedge\ldots\wedge d\phi_1\wedge d\psi_1.
	\end{align}
	For the rest of this proof, we will be calculating $\det((e_{ij}))$.	
	
	To be able to find the matrix $(e_{ij})$, we need to calculate Equation (\ref{eq:Proof_Haar_zo_ziet_e^k_er_in_praktijk_uit}) explicitly. Let $x\in G$ be a representative of $xK\in G/K$, then $x=\prod_{2\leq i\leq N}A(i)(\phi_{i-1},\psi_{i-1})$ for given $\phi_j,\psi_j$. If we write recursively $x_{n+1}(\phi_1,\ldots\phi_n,\psi_1,\ldots,\psi_n)=x_n(\phi_1,\ldots,\phi_{n-1},\psi_1,\ldots,\psi_{n-1})e^{\phi_{n}\lambda_{3}}e^{\psi_n\lambda_{n^2+1}}$ for $n\in\N$, we get $$x\equiv x_{N}=x_{N-1}e^{\phi_{N-1}\lambda_{3}}e^{\psi_{N-1}\lambda_{(N-1)^2+1}}.$$
	If we label $\omega_{x_l}:=x_{l}^{-1}dx_l$ for $1\leq l\leq N-1$, we get
	\begin{align}\label{eq:Proof_Haar_measure_omega_h_l}
	\omega_{x_{l+1}}&=e^{-\psi_{l}\lambda_{l^2+1}}e^{-\phi_{l}\lambda_{3}}x_{l}^{-1}dx_le^{\phi_l\lambda_{3}}e^{\psi_l\lambda_{l^2+1}}+e^{-\psi_l\lambda_{l^2+1}}d\phi_l\lambda_{3}e^{\psi_l\lambda_{l^2+1}}+d\psi_l\lambda_{l^2+1}\nonumber\\
	&=e^{-\psi_{l}\lambda_{l^2+1}}e^{-\phi_{l}\lambda_{3}}\omega_{x_l}e^{\phi_l\lambda_{3}}e^{\psi_l\lambda_{l^2+1}}+e^{-\psi_l\lambda_{l^2+1}}d\phi_l\lambda_{3}e^{\psi_l\lambda_{l^2+1}}+d\psi_l\lambda_{l^2+1}.
	\end{align}
	Putting Equation (\ref{eq:Proof_Haar_measure_omega_h_l}), with $l=N-1$, into Equation (\ref{eq:Proof_Haar_zo_ziet_e^k_er_in_praktijk_uit}) gives
	\begin{align*}
	(e^m)_{x1}&=\frac{1}{2}\Tr(\omega_{x_N} \lambda_{(N-1)^2+m-1})\\
	&=\frac{1}{2}\Tr\left(e^{-\psi_{N-1}\lambda_{(N-1)^2+1}}e^{-\phi_{N-1}\lambda_{3}}\omega_{x_{N-1}}e^{\phi_{N-1}\lambda_{3}}e^{\psi_{N-1}\lambda_{(N-1)^2+1}}\lambda_{(N-1)^2+m-1}\right)\\
	&\quad+\frac{1}{2}\Tr\left(e^{-\psi_{N-1}\lambda_{(N-1)^2+1}}d\phi_{N-1}\lambda_{3}e^{\psi_{N-1}\lambda_{(N-1)^2+1}}\lambda_{(N-1)^2+m-1}\right)+d\psi_{N-1}\,\delta_{(N-1)^2+1,(N-1)^2+m-1}.
	\end{align*}
	After a quick calculation, we see
	$$e^{-\psi_{N-1}\lambda_{(N-1)^2+1}}d\phi_{N-1}\lambda_{3}e^{\psi_{N-1}\lambda_{(N-1)^2+1}}=d\phi_{N-1}\sin(\psi_{N_1})\cos(\psi_{N-1})\lambda_{(N-1)^2}+O(\mathrm{diag}),$$ where by the notation $O(\mathrm{diag})$ we mean a diagonal matrix that can be disregarded when taking the trace form with $\lambda_{(N-1)^2+k-1}$. Therefore
	\begin{align*}
	(e^m)_{x1}&=d\psi_{N-1}\,\delta_{(N-1)^2+1,(N-1)^2+m-1}+\cos(\psi_{N-1})\sin(\psi_{N-1})d\phi_{N-1}\delta_{(N-1)^2,(N-1)^2+m-1}\\
	&\quad+\frac{1}{2}\Tr\left(e^{-\psi_{N-1}\lambda_{(N-1)^2+1}}e^{-\phi_{N-1}\lambda_{3}}\omega_{x_{N-1}}e^{\phi_{N-1}\lambda_{3}}e^{\psi_{N-1}\lambda_{(N-1)^2+1}}\lambda_{(N-1)^2+m-1}\right).
	\end{align*} 
	Remember that we were interested in $\det((e_{ij}))$ where $e^m=\sum_{j}e_{mj}dx_j$. The set $\{\psi_{N-1},\phi_{N-1},\ldots,\psi_1,\phi_1\}$ is a set of charts of $G/K$, so tracking down all the $d\psi_{N-1},d\phi_{N-1},\ldots,d\psi_1,d\phi_1$ gives the following matrix
	\begin{align*}
	(e_{ij})=\begin{pmatrix}
	1&0&\mathbf{0}^T\\
	0&\cos(\psi_{N-1})\sin(\psi_{N-1})&\frac{1}{2}\Tr\left(e^{-\phi_{N-1}\lambda_3}\omega_{x_{N-1}}e^{\phi_{N-1}\lambda_3}\lambda_{(N-1)^2+1}\right)\\
	\mathbf{0}&\mathbf{0}&\frac{1}{2}\Tr\left(e^{-\psi_{N-1}\lambda_{(N-1)^2+1}}e^{-\phi_{N-1}\lambda_3}\omega_{x_{N-1}}e^{\phi_{N-1}\lambda_3}e^{\psi_{N-1}\lambda_{(N-1)^2+1}}\Lambda\right)
	\end{pmatrix}.
	\end{align*}
	where $\Lambda$ is a vector with indices $\Lambda_{j-2}=\lambda_{(N-1)^2+j-1}$, $j=3,\ldots,2(N-1)$ and $\mathbf{0}$ is the $2(N-2)$ dimensional 0 vector. If $N=2$ we are done and see that $\det((e_{ij}))=\cos(\psi_{N-1})\sin(\psi_{N-1})$. If $N\geq 3$ we see that the lower right corner of $(e_{ij})$ is itself a $2(N-2)\times 2(N-2)$ matrix . Taking the determinant, and using $\Tr(AB)=\Tr(BA)$, is then
	\begin{align*}
	\det((e_{ij}))=\cos(\psi_{N-1})\sin(\psi_{N-1})\det\left(\frac{1}{2}\Tr\left(e^{-\phi_{N-1}\lambda_3}\omega_{x_{N-1}}e^{\phi_{N-1}\lambda_3}e^{\psi_{N-1}\lambda_{(N-1)^2+1}}\Lambda e^{-\psi_{N-1}\lambda_{(N-1)^2+1}}\right)\right).
	\end{align*}
	Calculating $e^{\psi_{N-1}\lambda_{(N-1)^2+1}}\lambda_{(N-1)^2+m-1}e^{-\psi_{N-1}\lambda_{(N-1)^2+1}}$ with $3\leq m\leq 2(N-1)$ we find the following relation 
	\begin{align*}
	e^{\psi_{N-1}\lambda_{(N-1)^2+1}}\lambda_{(N-1)^2+m-1}e^{-\psi_{N-1}\lambda_{(N-1)^2+1}} = (-1)^{m+1}\sin(\psi_{N-1})\lambda_{j(m)}+\cos(\psi_{N-1})\lambda_{(N-1)^2+m-1},
	\end{align*}
	where
	\begin{align*}
	j(m)=\begin{cases} \left(\frac{m-1}{2}\right)^2 \qquad\qquad\qquad\text{ if $m$ is odd,}\\
	\left(\frac{m-2}{2}\right)^2+1\qquad\;\,\qquad\text{ if $m$ is even.}
	\end{cases}
	\end{align*}
	Since the elements $\{\lambda_{j}\}_{j=1,\ldots,N^2-1}$ are orthogonal with respect to the trace, and $\omega_{x_{N-1}}$ has values in $\mathfrak{su}(N-1)$, we see that only the $\sin(\psi_{N-1})$ part contributes. Filling this in gives
	\begin{align}\label{eq:Proof_Haar_Measure_induction_omega_ h}
	\det((e_{ij}))&=\cos(\psi_{N-1})\sin^{2(N-1)-1}(\psi_{N-1})\det\left(\frac{1}{2}\Tr\left(e^{-\phi_{N-1}\lambda_3}\omega_{x_{N-1}}e^{\phi_{N-1}\lambda_3}\left((-1)^{m+1}\lambda_{j(m)}\right)_{m=3,\ldots,2(N-1)}\right)\right).
	\end{align}
	To finish the proof, we make the following claim:
	\begin{claim}
		$$\det\left(\frac{1}{2}\Tr\left(e^{-\phi_{N-1}\lambda_3}\omega_{x_{N-1}}e^{\phi_{N-1}\lambda_3}\left((-1)^{m+1}\lambda_{j(m)}\right)_{m=3,\ldots,2(N-1)}\right)\right) = 2 \prod_{j=1}^{N-2}\cos^{2j-1}(\psi_j)\sin(\psi_j).$$
	\end{claim}
	\begin{proof}
		
		Let $1<l\leq N-2$ and consider $$\frac{1}{2}\Tr\left(e^{-\phi_{l+1}\lambda_3}\omega_{x_{l+1}}e^{\phi_{l+1}\lambda_3}\lambda_{j(m)}\right)$$ for a given $m$. We apply Equation (\ref{eq:Proof_Haar_measure_omega_h_l}) and find 
		\begin{align}\label{eq:Proof_Haar_measure_recursive_equation}
		\frac{1}{2}\Tr\left(e^{-\phi_{l+1}\lambda_3}\omega_{x_{l+1}}e^{\phi_{l+1}\lambda_3}\lambda_{j(m)}\right) &=\frac{1}{2}\Tr\left(e^{-\phi_{l+1}\lambda_3}e^{-\psi_{l}\lambda_{l^2+1}}e^{-\phi_{l}\lambda_{3}}\omega_{x_{l}}e^{\phi_{l}\lambda_{3}}e^{\psi_{l}\lambda_{l^2+1}}e^{\phi_{l+1}\lambda_3}\lambda_{j(m)}\right)\nonumber\\ 
		&\quad+\frac{1}{2}\Tr\left(e^{-\phi_{l+1}\lambda_3}e^{-\psi_{l}\lambda_{l^2+1}}d\phi_{l}\lambda_{3}e^{\psi_{l}\lambda_{l^2+1}}e^{\phi_{l+1}\lambda_3}\lambda_{j(m)}\right)\nonumber\\
		&\quad+\frac{1}{2}\Tr\left(e^{-\phi_{l+1}\lambda_3}d\psi_{l}\lambda_{l^2+1}e^{\phi_{l+1}\lambda_3}\lambda_{j(m)}\right)\nonumber\\
		\begin{split}	
		&=\frac{1}{2}\Tr\left(e^{-\phi_{l}\lambda_{3}}\omega_{x_{l}}e^{\phi_{l}\lambda_{3}}\Ad(e^{\psi_{l}\lambda_{l^2+1}}e^{\phi_{l+1}\lambda_3})(\lambda_{j(m)})\right)\\ 
		&\quad+\frac{d\phi_{l}}{2}\Tr\left(\Ad(e^{-\phi_{l+1}\lambda_3}e^{-\psi_{l}\lambda_{l^2+1}})(\lambda_{3})\cdot\lambda_{j(m)}\right)\\
		&\quad+\frac{d\psi_{l}}{2}\Tr\left(\Ad(e^{-\phi_{l+1}\lambda_3})(\lambda_{l^2+1})\cdot\lambda_{j(m)}\right).
		\end{split}
		\end{align}
		To be able to evaluate this, we need some relations for the adjoint action. These are given by:
		\begin{align}\label{eq:recursion_eq_1}
		\Ad(e^{-\phi_n\lambda_3})\lambda_{q^2+1} =& \cos(\phi_n)\lambda_{q^2+1}-\sin(\phi_n)\lambda_{q^2},\\
		\Ad(e^{-\phi_n\lambda_3}e^{-\psi_{q}\lambda_{q^2+1}})\lambda_3=& \cos(\psi_{q})\sin(\psi_{q})\left[\cos(\phi_n)\lambda_{q^2}+\sin(\phi_n)\lambda_{q^2+1}\right]+O(\mathrm{diag})
		,\\
		\begin{split}
		\Ad\left(e^{\psi_{p}\lambda_{p^2+1}}e^{\phi_{n}\lambda_{3}}\right)\lambda_{q^2} =& \cos(\psi_p)\left[\cos(\phi_n)\lambda_{q^2}-\sin(\phi_n)\lambda_{q^2+1}\right]\\
		&-\sin(\psi_p)\left[\cos(\phi_n)\lambda_{p^2+2q}+\sin(\phi_n)\lambda_{p^2+2q+1}\right],
		\end{split}\\
		\begin{split}\label{eq:recursion_eq_4}
		\Ad\left(e^{\psi_{p}\lambda_{p^2+1}}e^{\phi_{n}\lambda_{3}}\right)\lambda_{q^2+1} =& \cos(\psi_p)\left[\sin(\phi_n)\lambda_{q^2}+\cos(\phi_n)\lambda_{q^2+1}\right]\\
		&+\sin(\psi_p)\left[-\sin(\phi_n)\lambda_{p^2+2q}+\cos(\phi_n)\lambda_{p^2+2q+1}\right],
		\end{split}
		\end{align}
		where $p,q,n\in \N$ and $p>q>1$. Before we fill this in, we recall that each $\lambda_j$ is orthogonal with respect to the trace form and $\omega_{x_l}$ has values in $\mathfrak{su}(l)$. Hence we see that only a few of the terms survive in Equations (\ref{eq:recursion_eq_1}-\ref{eq:recursion_eq_4}), and the only relevant terms are given here: 
		\begin{align}\label{eq:recursion_eq_1.2}
		\Ad(e^{-\phi_n\lambda_3})\lambda_{q^2+1} =& \cos(\phi_n)\lambda_{q^2+1}-\sin(\phi_n)\lambda_{q^2},\\
		\Ad(e^{-\phi_n\lambda_3}e^{-\psi_{q}\lambda_{q^2+1}})\lambda_3 =& \cos(\psi_{q})\sin(\psi_{q})\left[\cos(\phi_n)\lambda_{q^2}+\sin(\phi_n)\lambda_{q^2+1}\right]
		,\\
		\begin{split}
		\Ad\left(e^{\psi_{p}\lambda_{p^2+1}}e^{\phi_{n}\lambda_{3}}\right)\lambda_{q^2} =& \cos(\psi_p)\left[\cos(\phi_n)\lambda_{q^2}-\sin(\phi_n)\lambda_{q^2+1}\right],
		\end{split}\\
		\begin{split}\label{eq:recursion_eq_4.2}
		\Ad\left(e^{\psi_{p}\lambda_{p^2+1}}e^{\phi_{n}\lambda_{3}}\right)\lambda_{q^2+1} =& \cos(\psi_p)\left[\sin(\phi_n)\lambda_{q^2}+\cos(\phi_n)\lambda_{q^2+1}\right].
		\end{split}
		\end{align}	
		The latter two equations can we written even more compactly, namely
		\begin{align*}
		\begin{pmatrix}
		\Ad\left(e^{\psi_{p}\lambda_{p^2+1}}e^{\phi_{n}\lambda_{3}}\right)\lambda_{q^2}\\
		\Ad\left(e^{\psi_{p}\lambda_{p^2+1}}e^{\phi_{n}\lambda_{3}}\right)\lambda_{q^2+1}
		\end{pmatrix} = \cos(\psi_p)\begin{pmatrix}
		\cos(\phi_{n})&-\sin(\phi_n)\\
		\sin(\phi_n)&\cos(\phi_n)
		\end{pmatrix}\begin{pmatrix}
		\lambda_{q^2}\\
		\lambda_{q^2+1}
		\end{pmatrix}.
		\end{align*}
		Therefore, we see that due to linearity of the trace form that 
		\begin{align}
		\begin{pmatrix}
		\Tr\left(e^{-\phi_{l}\lambda_{3}}\omega_{x_{l}}e^{\phi_{l}\lambda_{3}}\Ad(e^{\psi_{l}\lambda_{l^2+1}}e^{\phi_{l+1}\lambda_3})\lambda_{q^2}\right)\\
		\Tr\left(e^{-\phi_{l}\lambda_{3}}\omega_{x_{l}}e^{\phi_{l}\lambda_{3}}\Ad(e^{\psi_{l}\lambda_{l^2+1}}e^{\phi_{l+1}\lambda_3})\lambda_{q^2+1}\right)
		\end{pmatrix} = \cos(\psi_l)\begin{pmatrix}
		\cos(\phi_{l+1})&-\sin(\phi_{l+1})\\
		\sin(\phi_{l+1})&\cos(\phi_{l+1})
		\end{pmatrix}\begin{pmatrix}
		\Tr\left(e^{-\phi_{l}\lambda_{3}}\omega_{x_{l}}e^{\phi_{l}\lambda_{3}}\lambda_{q^2}\right)\\
		\Tr\left(e^{-\phi_{l}\lambda_{3}}\omega_{x_{l}}e^{\phi_{l}\lambda_{3}}\lambda_{q^2+1}\right)
		\end{pmatrix}. 
		\end{align}
		Filling these equations in into Equation (\ref{eq:Proof_Haar_measure_recursive_equation}) gives
		\begin{align*}
		\frac{1}{2}\Tr\left(e^{-\phi_{l+1}\lambda_3}\omega_{x_{l+1}}e^{\phi_{l+1}\lambda_3}\lambda_{j(m)}\right) &=\frac{1}{2}\Tr\left[e^{-\phi_{l}\lambda_{3}}\omega_{x_{l}}e^{\phi_{l}\lambda_{3}}\Ad\left(e^{\psi_{l}\lambda_{l^2+1}}e^{\phi_{l+1}\lambda_3}\right)\lambda_{j(m)}\right]\\ 
		&\quad+d\phi_{l}\cos(\psi_{l})\sin(\psi_l)\left[\cos(\phi_{l+1})\delta_{l^2,j(m)}+\sin(\phi_{l+1})\delta_{l^2+1,j(m)}\right]\\
		&\quad+d\psi_{l}\left[\cos(\phi_{l+1})\delta_{l^2+1,j(m)}-\sin(\phi_{l+1})\delta_{l^2,j(m)}\right].
		\end{align*}
		Note that $j(m)$ is either a square number or a square number plus 1, which must mean that $\delta_{l^2,j(m)}$ can only be non-zero if $m$ is odd, and $\delta_{l^2+1,j(m)}$ can only be non-zero if $m$ is even.
		
		Define the $2l\times 2l$ dimensional matrix $$(X_{l+1})_{km}:=\frac{1}{2}\Tr\left(e^{-\phi_{l+1}\lambda_3}(\omega_{x_{l+1}})_k\;e^{\phi_{l+1}\lambda_3}\cdot(-1)^{m+3}\lambda_{j(m+2)}\right).$$ To prove the claim, we need to calculate $\det(X_{N-1})$. Swapping four rows and four colomns does not change the value of $\det(X_{N-1})$, so we swap the first and second row and colomn with the $(2(N-2)-1)$-th and the $2(N-2)$-th row and colomn respectively. Redefining this again as $X_{N-1}$, we get, using the above equations:
		
		\begin{align*}
		X_{N-1}&=\left(
		\scalemath{0.7}{\begin{array}{ccc}
			-\sin(\phi_{N-1}) & \cos(\phi_{N-1})\cos(\psi_{N-2})\sin(\psi_{N-2})&0\\
			-\cos(\phi_{N-1}) & -\sin(\phi_{N-1})\cos(\psi_{N-2})\sin(\psi_{N-2})&0\\
			0&0& \cos(\psi_{N-2})\begin{pmatrix}
			\cos(\phi_{N-1})&\sin(\phi_{N-1})\\
			-\sin(\phi_{N-1})&\cos(\phi_{N-1})
			\end{pmatrix}\begin{pmatrix}
			\frac{1}{2}\Tr\left(e^{-\phi_{N-2}\lambda_{3}}\omega_{x_{N-2}}e^{\phi_{N-2}\lambda_{3}}\lambda_{(\frac{m+1}{2})^2}\right)\\
			\frac{1}{2}\Tr\left(-e^{-\phi_{N-2}\lambda_{3}}\omega_{x_{N-2}}e^{\phi_{N-2}\lambda_{3}}\lambda_{(\frac{m}{2})^2+1}\right)
			\end{pmatrix}
			\end{array}}\right)\\
		&=\left(
		\scalemath{0.7}{\begin{array}{ccc}
			-\sin(\phi_{N-1}) & \cos(\phi_{N-1})\cos(\psi_{N-2})\sin(\psi_{N-2})&0\\
			-\cos(\phi_{N-1}) & -\sin(\phi_{N-1})\cos(\psi_{N-2})\sin(\psi_{N-2})&0\\
			0&0& \cos(\psi_{N-2})\begin{pmatrix}
			\cos(\phi_{N-1})&\sin(\phi_{N-1})\\
			-\sin(\phi_{N-1})&\cos(\phi_{N-1})
			\end{pmatrix}X_{N-2}
			\end{array}}\right).
		\end{align*}
		Taking the determinant gives 
		\begin{align*}
		\det(X_{N-1})=&\cos(\psi_{N-2})\sin(\psi_{N-2})\cos^{2(N-2)-2}(\psi_{N-2})\left[\det\begin{pmatrix}
		\cos(\phi_{N-1})&\sin(\phi_{N-1})\\
		-\sin(\phi_{N-1})&\cos(\phi_{N-1})
		\end{pmatrix}\right]^{N-3}\det\left(X_{N-2}\right)\\
		=&\cos^{2(N-2)-1}(\psi_{N-2})\sin(\psi_{N-2})\det\left(X_{N-2}\right).
		\end{align*}
		Recursively continuing the decomposition of the latter determinant gives
		\begin{align}\label{eq:Proof_Haar_measure_almost_there}\begin{split}
		\det\left(\frac{1}{2}\Tr\left(e^{-\phi_{N-1}\lambda_3}\right.\right.&\left.\omega_{x_{N-1}}e^{\phi_{N-1}\lambda_3}((-1)^{m+1}\lambda_{j(m)})_{m=3,\ldots,2(N-1)}\right)\bigg)\\
		&=\left[\prod_{j=2}^{N-2}\cos^{2j-1}(\psi_j)\sin(\psi_j)\right]\det\left(\frac{1}{2}\Tr(e^{-\phi_1\lambda_3}\omega_{x_2}e^{\phi_{1}\lambda_3}((-1)^{m+1}\lambda_{j(m)})_{m=3,4})\right).\end{split}
		\end{align}
		The last determinant is easily found, for $x_{2}=e^{\phi_1\lambda_3}e^{\psi_1\lambda_2}$ and so $$\omega_{x_2}=x_{2}^{-1}dx_{2}=d\phi_1 e^{-\psi_1\lambda_2}\lambda_3e^{\psi_1\lambda_2}+d\psi_1\lambda_2.$$ Since $\lambda_{j(3)}=\lambda_1$ and $\lambda_{j(4)}=\lambda_2$ we can find the final trace by just computing the matrix multiplications, which gives
		\begin{align*}
		\det\left(\frac{1}{2}\Tr(e^{-\phi_1\lambda_3}\omega_{x_2}e^{\phi_{1}\lambda_3}((-1)^{m+1}\lambda_{j(m)})_{m=3,4})\right) &= \det\begin{pmatrix}
		-\sin(2\phi_1)&\cos(2\phi_1)\sin(2\psi_1)\\
		-\cos(2\phi_1)&-\sin(2\phi_1)\sin(2\psi_1)
		\end{pmatrix}\\
		&=\sin(2\psi_1)=2\sin(\psi_1)\cos(\psi_1).
		\end{align*}
		Filling this in into Equation (\ref{eq:Proof_Haar_measure_almost_there}) gives the result	
		\begin{align}\label{eq:Proof_Haar_measure_claim_equation}
		\det\left(\frac{1}{2}\Tr\left(e^{-\phi_{N-1}\lambda_3}\omega_{x_{N-1}}e^{\phi_{N-1}\lambda_3}\left((-1)^{m+1}\lambda_{j(m)}\right)_{m=3,\ldots,2(N-1)}\right)\right) = 2\prod_{j=1}^{N-2}\cos^{2j-1}(\psi_j)\sin(\psi_j)
		\end{align} which proves the claim.
	\end{proof}
	Putting Equation (\ref{eq:Proof_Haar_measure_claim_equation}) and (\ref{eq:Proof_Haar_Measure_induction_omega_ h}) into Equation (\ref{eq:Proof_Haar_measure_equation_that_only_matters}) gives
	\begin{align*}
	dg&=dg_{K}dk\\
	&= \det((e_{ij}))\;d\phi_1\ldots d\phi_{N-1}d\psi_1\ldots d\psi_{N-1} dg_{SU(N-1)} d\omega_N\\
	&=2\cos(\psi_{N-1})\sin^{2(N-1)-1}(\psi_{N-1})\left[\prod_{j=1}^{N-2}\cos^{2j-1}(\psi_{j})\sin(\psi_j)\right]d\phi_1\ldots d\phi_{N-1}d\psi_1\ldots d\psi_{N-1} dg_{SU(N-1)}d\omega_{N-1}.
	\end{align*}
	Thus this is the Haar measure up to a normalization constant. To get the normalised Haar measure, we need to explicitely integrate over the whole group. The normalisation constant $C_N$ in Equation (\ref{eq:Haar_Measure}) can be found by noting that the only non-trivial integration is over the $\psi_j$ coordinates, and each integral can be evaluated using the following identity $$\int_0^{\frac{\pi}{2}}\sin^k(x)\cos(x)dx=\int_0^{\frac{\pi}{2}}\cos^k(x)\sin(x)dx=\frac{1}{k+1}.$$

\end{proof}

\nocite{*}

\newpage
\bibliography{bibfile}
	
\end{document}